\documentclass[10pt]{amsart}
\usepackage{color}
\usepackage{verbatim}
\usepackage[pdftex]{graphicx}

\usepackage{amssymb}
\allowdisplaybreaks[4]

\def\H{\mathcal{H}}

\def\LT{\left}
\def\RT{\right}
\def\oo{(1+o(1))}

\definecolor{c20}{rgb}{0.,0.7,0.}
\definecolor{c30}{rgb}{0.,0.,1.}
\definecolor{c40}{rgb}{1,0.1,0.7}
\definecolor{c50}{rgb}{1,0,0}
\definecolor{c60}{rgb}{1,0.9,0.1}

\def\x{\vk{x}}

\def\SI{\Sigma}

\newcommand{\abs}[1]{\left\lvert #1 \right\rvert}

\newcommand{\E}[1]{\mathbb{E}\left\{ #1\right\}}

\newcommand{\pk}[1]{\mathbb{P} \left\{ #1 \right \} }

\newcommand{\R}{\mathbb{R}}

\newcommand{\inr}{\in \R}

\newcommand{\ldot}{,\ldots,}

\newcommand{\limit}[1]{\lim_{#1 \to   \infty}}

\newcommand{\BQN}{\begin{eqnarray}}
\newcommand{\EQN}{\end{eqnarray}}
\newcommand{\BQNY}{\begin{eqnarray*}}
\newcommand{\EQNY}{\end{eqnarray*}}

\newcommand{\BS}{\begin{sat}}
\newcommand{\ES}{\end{sat}}
\newcommand{\BT}{\begin{theo}}
\newcommand{\ET}{\end{theo}}
\newcommand{\BK}{\begin{korr}}
\newcommand{\EK}{\end{korr}}

\newcommand{\BD}{\begin{de}}
\newcommand{\ED}{\end{de}}

\newcommand{\BIT}{\begin{itemize}}
\newcommand{\EIT}{\end{itemize}}
\newcommand{\BDI}{\begin{description}}
\newcommand{\EDI}{\end{description}}

\newcommand{\BRM}{\begin{remarks}}
\newcommand{\ERM}{\end{remarks}}

\newcommand{\BEL}{\begin{lem}}
\newcommand{\EEL}{\end{lem}}

\newtheorem{theo}{Theorem}[section]
\newtheorem{sat}[theo]{Proposition}
\newtheorem{de}[theo]{Definition}
\newtheorem{lem}[theo]{Lemma}

\newtheorem{korr}[theo]{Corollary}

\newtheorem{remarks}[theo]{Remarks}
\newtheorem{prop}[theo]{Proposition}

\newcommand{\nelem}[1]{{Lemma \ref{#1}}}

\newcommand{\netheo}[1]{{Theorem \ref{#1}}}

\newcommand{\COM}[1]{}

\newcommand{\QED}{\hfill $\Box$}

\newcommand{\kb}[1]{\boldsymbol{#1}}
\newcommand{\vk}[1]{\kb{#1}}

\def\IF{\infty}

\def\LT{\left}
\def\RT{\right}

\def\vn{\varepsilon}

\def\Del{\triangle}

 \def\squ {\sqrt{u}}

\def\K1#1{\textcolor{cyan}{#1}}

\def\K1#1{\textcolor{cyan}{#1}}

\def\LLc#1{\textcolor{black}{#1}}
\def\kkk#1{\textcolor{black}{#1}}

\def\jj#1{\textcolor{black}{#1}}
\def\kkkk#1{\textcolor{black}{#1}}

\def\oA{\overline A}
\def\oB{\overline B}
\def\mcE{U}

\begin{document}

\title{Exact asymptotics of component-wise extrema of two-dimensional Brownian motion }
\author{Krzysztof D\c{e}bicki}
\address{Krzysztof D\c{e}bicki, Mathematical Institute, University of Wroc\l aw, pl. Grunwaldzki 2/4, 50-384 Wroc\l aw, Poland}
\email{Krzysztof.Debicki@math.uni.wroc.pl}

\author{Lanpeng Ji}
\address{Lanpeng Ji, School of Mathematics, University of Leeds, Woodhouse Lane, Leeds LS2 9JT, United Kingdom
}
\email{l.ji@leeds.ac.uk}

\author{Tomasz Rolski}
\address{Tomasz Rolski, Mathematical Institute, University of Wroc\l aw, pl. Grunwaldzki 2/4, 50-384 Wroc\l aw, Poland}
\email{Tomasz.Rolski@math.uni.wroc.pl}


\bigskip

\date{\today}
 \maketitle

 {\bf Abstract:}
We derive the exact asymptotics of
\[
\pk{  \sup_{t\ge 0} \Bigl( X_1(t) -  \mu_1 t\Bigr)>  u, \  \sup_{s\ge 0} \Bigl( X_2(s) -  \mu_2 s\Bigr)>  u  },\ \   u\to\infty,
\]
where
$(X_1(t),X_2(s))_{t,s\ge0}$
is a correlated two-dimensional Brownian motion with correlation $\rho\in[-1,1]$
and $\mu_1,\mu_2>0$.
It appears that the play between $\rho$ and $\mu_1,\mu_2$
leads to several types of asymptotics.
Although  the exponent in the asymptotics as a function of $\rho$ is continuous,
\jj{one} can observe different types of prefactor functions depending on the range of $\rho$,
which constitute a  phase-type transition phenomena.

 {\bf Key Words:} Two-dimensional Brownian motion;  exact asymptotics;  \jj{component-wise extrema}; quadratic programming problem;
 generalised Pickands-Piterbarg constants.

 {\bf AMS Classification:} Primary 60G15; secondary 60G70

 \section{Introduction}

Distributional properties of component-wise extrema of stochastic processes
attract growing interest in recent literature.
On one side, it is a natural object of interest in the extreme value theory of random fields.
On the other side, strong motivation to investigate  component-wise extrema stems for example from
multivariate stochastic models applied to modern multidimensional risk theory, financial mathematics
or advanced communication networks, to name some of the applied-probability areas.

We consider a standard correlated
 Brownian motion $(X_1(t),X_2(t))_{t\ge0}$ with constant correlation
   $\rho\in[-1,1]$, and let $(X_1(t),X_2(s))_{t,s\ge0}$ be its  two parameter extension, where
   $$\E{ X_1(t)X_2(s)}=\rho \min(t,s).$$
The aim of this paper is to find exact asymptotics of
   \begin{equation}\label{eq:PPu}
  P(u):=   \pk{Q_1>u,Q_2>u},\qquad u\to\infty,
   \end{equation}
   where $Q_j=\sup_{t\ge0}(X_j(t)-\mu_jt)$ with  $\mu_j>0$, $j=1,2$.

Due to its importance in, e.g., quantitative finance or ruin theory, 
the component-wise maxima 
$$(Q_1(T), Q_2(T))=
   \LT(\sup_{t\in[0,T]} (X_1(t)-\mu_1 t),  \sup_{s\in[0,T]} (X_2(s)-\mu_2 s)\RT)$$
have been studied extensively;
see, e.g.,  \cite{KZ16, HKR98, Met10, ShaoWang13, RS06}.
In particular, some formulas for the joint distribution of 
$(Q_1(T),\ Q_2(T))$
are known.
Unfortunately, they are in the form of infinite-sums of integrals of some special
functions, \jj{which makes them of limited use  in drawing out  asymptotic properties of
$P(u)$ as $u\to\infty$.}

Interestingly, in  \cite{KZ16} it was worked out a formula for joint survival function of
$(Q_1(\mathcal E_p), Q_2(\mathcal E_p))$,
where $\mathcal E_p$ is an independent   exponential random variable with parameter $p>0$.
Vector $(Q_1(\mathcal E_p), Q_2(\mathcal E_p))$ as well as $(Q_1,Q_2)$ have bivariate exponential distribution (BVE)
in the sense of the terminology of Kou and Zhong \cite{KZ16}, that is:
(i) it has exponential marginals and
(ii) it is absolute continuous with respect to two-dimensional Lebesgue measure.
The later property \jj{for $(Q_1,Q_2)$} follows from
Theorem 7.1 in \cite{AzW09} combined with the fact that $\pk{Q_j=0}=0$, see also related Lemma 4.4 in \cite{DHW19}.
We remark that requirement (ii) implies that $(Q_1,Q_2)$
does not belong to the  classical examples of  {\it Marshall-Olkin}-type  BVE; see \cite{MO67}.
Since there are no results in the literature on qualitative properties of our BVE distribution,
as a by-product of the \jj{results of this contribution}, we  analyze the
dependence structure of $Q_1$ and $Q_2$
in an asymptotical sense of Resnick \cite{Res87}; see Remarks \ref{Rem:main1} (b) and Remarks \ref{Rem:main2} (b)
for more details.
We refer also to a related work of Rogers and Shepp  \cite{RS06} who considered correlation structure of $(Q_1(T),Q_2(T))$
for two Brownian motions without drift.

 A need to consider the joint survival function for $(Q_1,Q_2)$ appeared
also in Lieshout and Mandjes \cite{LM07} who considered two parallel queues sharing the same Brownian input
(which is the case of $\rho=1$) and also a Brownian  tandem queue.
We refer to \cite{Mandjes} for further discussions on Gaussian-related queueing models
and
to \cite{DHJT18, DHM19}  for the analysis of a related
{\it simultaneous ruin} problem for the correlated Brownian motion model.

 It is worth noting that in recent papers \cite{HH19, Honnappa}, the component-wise maxima in discrete models defined by
$$ (\max_{1\le i \le n}X_i^1,\ldots, \max_{1\le i \le n}X_i^d),$$
with $(X_i^1,\ldots,X^d_i)$ $(i=1,2,\ldots)$  independent and identically distributed Gaussian random vectors, were discussed.

\smallskip

The first step in understanding \jj{the} asymptotics of \eqref{eq:PPu} is to find its logarithmic asymptotics. This was done recently in
\cite{DJT19a}, in an insurance  context, where $P(u)$  was interpreted as the probability of \jj{component-wise ruin.} 
 More precisely, by an application of Theorem 1 in \cite{RolskiSPA} 
\BQN\label{eq:logri}
\frac{\ln P(u)}{u} \sim -\frac{g(\vk{t}_0)}{2}, \ \ \ u\to\IF,
\EQN
where
\BQN\label{eq:twolayer}
g(\vk{t}_0)=\inf_{\vk t> \vk 0} \inf_{\vk{v} \ge \vk{1} + \vk{\mu} \vk t}  \vk{v} ^\top \Sigma_{ts}^{-1}   \vk{v}
\EQN
and $\Sigma_{ts}^{-1}$ is the inverse matrix of
$
\  \Sigma_{ts}=
\left(
   \begin{array}{cc}
     t & \rho\  t\wedge s \\
     \rho\  t\wedge s & s \\
   \end{array}
 \right),
$
with $\vk t=(t,s)^\top$ and $ t\wedge s=\min(t,s)$.  The main contribution of \cite{DJT19a} includes
the detailed analysis of the two-layer minimisation problem involved in $g(\vk t_0)$,
which results in an explicit logarithmic asymptotics of $P(u)$; see also Proposition \ref{Lem:Optg} below.

In order to get  the exact asymptotics of $P(u)$ as $u\to\infty$,
we employ a modification of the {\it double-sum} technique,
accommodated to the analysis of multivariate extremes investigated  in this contribution;
see Theorems \ref{Thm:mu12} and \ref{Thm:mu}, which constitute the main results of this paper.
It appears that the play between $\rho$ and $\mu_1,\mu_2$
leads to several types of  asymptotics. Although in  \cite{DJT19a} it was noticed,
that the exponent in the asymptotics as a function of $\rho$, called \jj{therein} an adjustment coefficient,
is continuous, \jj{one} can observe different types of prefactor  functions depending on the range of $\rho$.
This  {\it phase}-type phenomena
has no intuitive explanations.

In the rest of the paper we assume that  $\rho\in(-1,1)$ and
without loss of generality suppose
that $\mu_1\le\mu_2$.
Note that for $\rho=1$,
\BQNY
P(u)=\pk{ \sup_{s\ge 0} (X_2(s)-\mu_2 s)>u}=e^{-2 \mu_2 u},\ \ \ \forall u>0
\EQNY
 and, for $\rho=0$,
\BQN\label{eq:Prho0}
P(u) = \pk{\sup_{t\ge 0} (X_1(t)-\mu_1 t)>u}\pk{ \sup_{s\ge 0} (X_2(s)-\mu_2 s)>u}
 = e^{-2(\mu_1+\mu_2)u},\ \ \ \forall u>0.
\EQN
To work out the case $\rho=-1$, one  can use  a result from
\cite{TeuGoo1994}, to show that 
 \BQN\label{eq:TG}
 P(u) \sim e^{-(2\mu_2+6\mu_1) u} (2 I_{(\mu_1=\mu_2)}+ I_{(\mu_1<\mu_2)}),\ \ \ u\to \IF,
 \EQN
where $I_{(\cdot)}$ is the indicator function.


The rest of this paper is organised as follows. In
Section \ref{Sec:main},
we present the  exact asymptotics of $P(u)$, given in  Theorems \ref{Thm:mu12}, \ref{Thm:mu}.
Section \ref{s.log} recalls the {explicit} expressions for $g(\vk{t}_0)$ and $\vk t_0$ derived in \cite{DJT19a}.
The main lines of proofs  are displayed in  Section \ref{proofmain1} and Section \ref{proofmain2}, respectively,
followed by the Appendix consisting of technical  calculations.

We conclude this section by showing some notation and conventions used in this work.
All vectors here are $2$-dimensional column vectors written in bold letters. For instance
$\vk{\alpha}=(\alpha_1,\alpha_2)^\top$, with $^\top$ the transpose sign. Operations with vectors are meant component-wise, so
$\lambda \x=\x \lambda = ( \lambda x_1, \lambda x_2)^\top$ for any $\lambda\inr, \x\inr^2$. %
For any set $D \subseteq  [0,\IF)^2$, any $\lambda>0$ and any $(a_1,a_2)\in  [0,\IF)^2$ denote
\BQNY
\lambda D =\{(\lambda t, \lambda s): (t,s)\in D\},\ \ \ {(a_1,a_2)}+D=\{(a_1+ t, a_2+ s): (t,s)\in D\}.
\EQNY

Next, let us briefly mention the following standard notation
for two given positive functions $f(\cdot)$ and $h(\cdot)$. We write
$ f(x)=h(x)(1+o(1))$  or simply $f(x)\sim h(x)$, if  $ \lim_{x \to a}  {f(x)}/{h(x)} = 1$ ($a\in\R\cup\{\IF\}$). Further, write $ f(x) = o(h(x))$ if $ \lim_{x \to a}  {f(x)}/{h(x)} = 0$, and  
write $f(x) \lesssim h(x) $ if $ \lim_{x \to a}  {f(x)}/{h(x)} \le 1.$

\section{Main results}\label{Sec:main}

{In this section we  present} the exact asymptotics of $P(u)$, for which we need some additional notation.
First, define
\BQN\label{eq:rho12}
\hat\rho_1=\frac{\mu_1+\mu_2-\sqrt{(\mu_1+\mu_2)^2-4\mu_1(\mu_2-\mu_1)}}{4\mu_1}\in [0,\frac{1}{2}), \ \ \ \ \
\hat\rho_2=\frac{\mu_1+\mu_2}{2\mu_2}.
\EQN
These are key points, based on which we consider different scenarios of $\rho$.
\COM{
\begin{prop}\label{p.log}
\begin{eqnarray*}
\frac{\log(P(u))}{u}\sim
\left\{
  \begin{array}{ll}
    2(\mu_2+(1-2\rho)\mu_1), & \hbox{if \ $-1< \rho  \le\hat \rho_1$;} \\
    \frac{\mu_1+\mu_2+2/t^*}{1+\rho}, & \hbox{if \ $\hat \rho_1  < \rho<\hat \rho_2$ ;} \\
    2\mu_2, & \hbox{if \ $\hat \rho_2\le\rho<1$.}
  \end{array}
\right.
\end{eqnarray*}
\end{prop}
The proof of Proposition \ref{p.log} is postponed to Section \ref{s.log}.

}
Next, let
\BQN\label{eq:Sig-b}
\  \Sigma_{*}=
\left(
   \begin{array}{cc}
     t^*  & \rho s^* \\
     \rho s^*  & s^*  \\
   \end{array}
 \right),\ \ \ \ \ \vk b_*=(1+\mu_1 t^*, 1+\mu_2 s^*)^\top,
\EQN
with
\BQN\label{eq:tsstar}
 t^*=t^*(\rho)= s^*= s^*(\rho):=\sqrt {\frac{2(1-\rho)}{\mu_1^2+\mu_2^2-2\rho\mu_1\mu_2}}.
\EQN
Moreover,  denote, for any fixed $T,S>0$, 
\BQN\label{eq:DelTS}
\Del_{T,S}=\{(t,s): t\in [0,T], s\in [t,t+S]\} \cup \{ (t,s): s\in [0,T], t\in [s,s+S]\},
\EQN
and
define
\BQNY
\H(T,S):=\int_{\R^2}e^{\vk x^\top \Sigma_{*}^{-1}  \vk{b}_{*} }\pk{  {\underset{ (t,s)\in \Del_{T,S}}\exists }
\begin{array}{ccc}
X_1(t)- \mu_1  t >  x_1 \\
X_2(s)-  \mu_2  s  >x_2
\end{array}
 } dx_1 dx_2\in(0,\IF),
\EQNY
where the finiteness can be proved by following a standard argument in proving the finiteness of Pickands and Piterbarg type constants; see, e.g.,
\cite{Pit96} (or Lemma 4.2 in \cite{DHJT18}). Interestingly, a new Pickands-Piterbarg constant
\BQNY
\widetilde\H:=\lim_{S\to\IF}\lim_{T\to\IF} \frac{1}{T}\H(T,S)\in(0,\IF)
\EQNY
appears in the scenario  $ \hat \rho_1  < \rho< \hat \rho_2$; the existence, finiteness and positiveness of this constant are proved in \netheo{Thm:mu12} below.

We split the statement of the main \jj{results} on the exact asymptotics into two scenarios:
 $\mu_1<\mu_2$ and $\mu_1=\mu_2$  respectively.

\BT \label{Thm:mu12} Suppose that  $\mu_1<\mu_2$. We have, as $u\to\IF,$
\begin{eqnarray*}
P(u)
\sim
\left\{
  \begin{array}{ll}
    e^{- 2(\mu_2+(1-2\rho)\mu_1)u}, & \hbox{if \ $-1< \rho  <\hat \rho_1$;} \\
    \frac{1}{2}\       e^{- 2(\mu_2+(1-2\hat\rho_1)\mu_1)u}, & \hbox{if \ $\rho = \hat \rho_1$;} \\
    \frac{\widetilde\H \sqrt{t^*}}{2\sqrt{\pi(1-\rho)}}\ u^{-1/2} e^{- \frac{\mu_1+\mu_2+2/t^*}{1+\rho} u}, & \hbox{if \ $ \hat \rho_1  < \rho< \hat \rho_2 $;}\\
     e^{- 2\mu_2 u}, & \hbox{if \ $ \hat \rho_2<\rho <1$,}
  \end{array}
\right.
\end{eqnarray*}
where
$$
0<\frac{t^* \vk\mu^\top \Sigma_*^{-1} \vk b_*}{16 \prod_{i=1}^2(\Sigma_*^{-1} \vk b_*)_i}<\widetilde\H<\IF.
$$



\ET
\begin{remarks} \label{Rem:main1}
(a). It turns out that the special scenario $\rho=\hat\rho_2$ is of different nature than the scenarios analyzed
in Theorem \ref{Thm:mu12}.
Note that in this case we have $b_1=b_2=0$ in Lemma \ref{Lem:Taylor}, which implies that
around \jj{its} optimizing  point $(t^*,s^*)=(1/\mu_2,1/\mu_2)$ function $g(t,s)$ \jj{defined in Section 3} takes different form
than for other scenarios. This makes its analysis go out of the approach that works for the other scenarios.
In Section \ref{sc.4}, following the same lines of reasoning as given in the proof of case $ \hat \rho_2<\rho <1$ in Theorem \ref{Thm:mu12},
 we find the following bounds for \jj{the case of} $\rho=\hat\rho_2$
 \BQN
 \frac{ 1}{2}\ e^{- 2\mu_2 u}\  \lesssim\  P(u)\ \lesssim \   e^{- 2\mu_2 u}, \ {\rm as}\ u\to\infty.\label{old.4}
 \EQN

(b). It follows from 
Theorem \ref{Thm:mu12} and \eqref{eq:TG} that for any $-1 \le\rho<\hat \rho_2$
\BQNY
\pk{Q_1(\IF)>u| Q_2(\IF)>u}=\frac{\pk{Q_1(\IF)>u, Q_2(\IF)>u}}{\pk{Q_2(\IF)>u}}\to 0,\ \ u\to\IF.
\EQNY
{ According to  the terminology from} \cite{Res87}, this means that
$Q_1(\IF)$ is  asymptotically independent of $Q_2(\IF)$.
Similarly, one can see that $Q_2(\IF)$ is also asymptotically independent of $Q_1(\IF)$
(note that the notion of asymptotically independence is not symmetric).
Furthermore,  for $  \hat \rho_2 \le \rho \le 1$, we have that  $Q_2(\IF)$ is asymptotically independent of
$Q_1(\IF)$, but $Q_1(\IF)$ is asymptotically dependent of (equivalent to) $Q_2(\IF)$.

\end{remarks}

Next we give the result for the case where $\mu:=\mu_1=\mu_2$. In this case,
we have $t^*=s^*=1/\mu$ and
\BQNY
\widetilde\H:=\lim_{S\to\IF}\lim_{T\to\IF} \frac{1}{T}\int_{\R^2}e^{\frac{2\mu}{1+\rho}(x_1+x_2)}\pk{  {\underset{ (t,s)\in \Del_{T,S}}\exists }
\begin{array}{ccc}
X_1(t)- \mu  t >  x_1 \\
X_2(s)-  \mu  s  >x_2
\end{array}
 } dx_1 dx_2.
\EQNY

\BT \label{Thm:mu} Suppose that $\mu_1=\mu_2$. We have, as $u\to\IF,$
\begin{eqnarray*}
P(u)
\sim
\left\{
  \begin{array}{ll}
    2 \ e^{- 4(1-\rho)\mu  u}, & \hbox{if \ $-1< \rho  < 0$;} \\
       \ e^{- 4\mu  u}, & \hbox{if \ $\rho = 0$;} \\
    \frac{\widetilde\H }{2\sqrt{\pi\mu(1-\rho)}}\ u^{-1/2} e^{- \frac{4\mu}{1+\rho}  u}, & \hbox{if \ $ 0  < \rho< 1$, }
      \end{array}
\right.
\end{eqnarray*}
where
$
(1+\rho)/{16 }<\widetilde \H<\IF.
$
\ET

\begin{remarks}\label{Rem:main2}
(a). Note that comparing scenario $-1< \rho  < 0$ of Theorem \ref{Thm:mu}
with $-1< \rho  <\hat \rho_1$ of Theorem \ref{Thm:mu12},  there is an additional 2 appearing in the
asymptotics. The reason for this is that there are two equally important minimizers of $g(t,s), (t,s)\in(0,\IF)^2$
in \LLc{the case of $\mu_1=\mu_2$.}


(b). For any $-1< \rho<1$, we have that $Q_1(\IF)$  and $Q_2(\IF)$ are mutually asymptotically independent.

\end{remarks}


\section{Analysis of the two-layer minimization problem }\label{s.log}

In this section, for completeness and for reference we recall some notation and the result on the two-layer minimization problem \LLc{\eqref{eq:twolayer}} derived in \cite{DJT19a}. Recall that
$g(\vk{t}_0)=\inf_{(t,s)\in(0,\IF)^2} g(t,s)
$   
with
$$g(t,s):= \underset{y\ge 1+\mu_2 s}{\inf_{x\ge 1+\mu_1 t}}\  (x,y)\ \Sigma_{ts}^{-1}\  (x,y)^\top,\ \ \ \ \ t,s>0.
$$
We define for $t,s>0$ the following functions:
\BQNY
&&g_1(t)=\frac{(1+\mu_1 t)^2}{t}, \ \ \ \ g_2(s)=\frac{(1+\mu_2 s)^2}{s},\\
 && g_3(t,s)=(1+\mu_1 t, 1+\mu_2 s) \ \Sigma_{ts}^{-1} \  (1+\mu_1 t, 1+\mu_2 s)^\top.
\EQNY
Since $t \wedge s$ appears in the above formula, we shall consider a partition of the quadrant $(0,\IF)^2$, namely
\BQN\label{eq:AB}
(0,\IF)^2=A\cup L\cup B,\ \ \ A=\{s< t\},\ L=\{s=t\},\  B=\{s>t\}.
\EQN
For convenience we denote $ \oA=\{s\le t\}=A\cup L$ and  $ \oB=\{s\ge t\}=B\cup L$. Hereafter, all   sets are defined on $(0,\IF)^2$, so $(t,s)\in (0,\IF)^2$ will be omitted.



Note that $g_3(t,s)$ can be represented in the following two different  forms:
\BQN
g_3(t,s)&=&\left\{\begin{array}{cc}
g_A(t,s):= \frac{ (1+\mu_1 t)^2 s -2\rho s  (1+\mu_1 t)(1+\mu_2 s) + (1+\mu_2s)^2 t }{ts-\rho^2 s^2}, & \hbox{if }  (t,s)\in \oA \label{eq:g31}\\
g_B(t,s):=\frac{ (1+\mu_1 t)^2 s -2\rho t  (1+\mu_1 t)(1+\mu_2 s) + (1+\mu_2s)^2 t }{ts-\rho^2 t^2},  & \hbox{if }   (t,s)\in \oB
\end{array}\right.\\
&=&\left\{\begin{array}{cc}
\frac{(1+\mu_2s)^2}{s} +\frac{((1+\mu_1 t)-\rho(1+\mu_2 s))^2}{t-\rho^2 s}, & \hbox{if }  (t,s)\in \oA \\
\frac{(1+\mu_1 t)^2}{t} +\frac{((1+\mu_2 s)-\rho(1+\mu_1 t))^2}{s-\rho^2 t}, & \hbox{if }   (t,s)\in \oB.
\end{array}\right.\label{eq:g32}
\EQN
Denote further
\BQN\label{eq:fL}
g_L(s):=g_A(s,s)=g_B(s,s)=\frac{(1+\mu_1s)^2 +(1+\mu_2s)^2-2\rho(1+\mu_1s)(1+\mu_2s) }{(1-\rho^2) s},\ \ \ s>0.
\EQN

The following result gives a full analysis of the two-layer minimization problem \LLc{\eqref{eq:twolayer}}, which is crucial for our derivation of the exact asymptotics of $P(u)$. 
We refer to \cite{DJT19a} for its detailed proof.

\begin{prop} \label{Lem:Optg} 
\begin{itemize}
\item[(i).] Suppose that  $-1< \rho < 0$. \\
For    $\mu_1<\mu_2$ we have
\BQNY
 g(\vk{t}_0) = g_A(t_A,s_A) =   4(\mu_2+(1-2\rho)\mu_1),
\EQNY
where,
$
(t_A,s_A)=(t_A(\rho),s_A(\rho)):=\LT(\frac{1-2\rho}{\mu_1}, \frac{1}{\mu_2-2\mu_1\rho} \RT)
\in A$ 
is the unique minimizer of $g(t,s), (t,s)\in(0,\IF)^2$.
\\For $\mu_1=\mu_2=:\mu$ we have
\BQNY
 g(\vk{t}_0) = g_A(t_A,s_A) = g_B(t_B,s_B)=  8(1-\rho)\mu,
\EQNY
where
$
(t_A,s_A)=\LT(\frac{1-2\rho}{\mu}, \frac{1}{ (1-2 \rho)\mu} \RT),
 (t_B,s_B):=\LT(\frac{1}{ (1-2 \rho)\mu},\frac{1-2\rho}{\mu} \RT)\in B$
are the only two  minimizers of $g(t,s), (t,s)\in(0,\IF)^2$.

\item[(ii).] Suppose that $0\le \rho<  \hat \rho_1$. We have
\BQNY
g(\vk{t}_0)  = g_A(t_A,s_A) =   4(\mu_2+(1-2\rho)\mu_1),
\EQNY
where $(t_A,s_A)$ is the unique minimizer  of $g(t,s), (t,s)\in(0,\IF)^2$.  

\item[(iii). ] Suppose that $\rho = \hat \rho_1$. We have
\BQNY
g(\vk{t}_0) = g_A(t_A,s_A) =   4(\mu_2+(1-2\rho)\mu_1),
\EQNY
where $(t_A,s_A) =(t_A(\hat\rho_1), s_A(\hat\rho_1))=(t^*(\hat\rho_1), s^*(\hat\rho_1))\in L$,
is the unique minimizer  of $g(t,s), (t,s)\in(0,\IF)^2$, with $(t^*,s^*)$ defined in \eqref{eq:tsstar}.

\item[(iv).] Suppose that $ \hat \rho_1  < \rho< \hat \rho_2 $. We have
\BQNY
g(\vk{t}_0)  = g_A(t^*, s^{*}) = g_L(t^*)= \frac{2}{1+\rho}(\mu_1+\mu_2+2/t^*),
\EQNY
where $(t^*,s^*)\in L$ is \LLc{the unique minimizer  of $g(t,s), (t,s)\in(0,\IF)^2$.}

\item[(v).] Suppose that $ \rho  =  \hat \rho_2 $. We have    $t^*(\hat\rho_2)=s^*(\hat\rho_2)=1/\mu_2$, and
\BQNY
g(\vk{t}_0) = g_A(1/\mu_2 ,1/\mu_2 )=g_L( 1/\mu_2 ) =g_2( 1/\mu_2 )=4\mu_2,
\EQNY
where
the  minimum  of $g(t,s), (t,s)\in(0,\IF)^2$ is attained at $(1/\mu_2 ,1/\mu_2 )$,  with
$g_3(1/\mu_2 ,1/\mu_2 )=g_2(1/\mu_2)$, and $1/\mu_2 $ is the unique minimizer of $g_2(s), s\in(0,\IF)$.

\item[(vi).] Suppose that  $ \hat \rho_2<\rho <1$. We have
\BQNY
g(\vk{t}_0)  = \inf_{(t,s)\in D_2}g_2(s)  =  g_2( {1}/{\mu_2})= 4\mu_2,
\EQNY
\kkk{where
the minimum  of $g(t,s), (t,s)\in(0,\IF)^2$ is attained when $g(t,s)=g_2(s)$. 
}
\end{itemize}
\end{prop}

%
%


\section{Proof of  Theorem \ref{Thm:mu12}} \label{proofmain1} 


The proof of Theorem \ref{Thm:mu12} will be presented in the order of cases
(i) $-1< \rho  <\hat \rho_1$, (ii) $ \hat \rho_1  < \rho< \hat \rho_2 $, (iii) $\rho = \hat \rho_1$, (iv) $ \hat \rho_2\le\rho <1$
in the following subsections.

Note that by self-similarity
\BQN \label{eq:SelfSimi}
P(u)=\pk{\exists_{t ,s> 0} \ (X_1(t)>(1+\mu_1 t)\sqrt u, \  X_2(s)>(1+\mu_2 s)\sqrt u)},
\EQN
and recall 
the notation for the optimizer points $(t_A,s_A)$  as introduced in Proposition \ref{Lem:Optg}. 

\subsection{(i) Scenario  $-1< \rho<\hat\rho_1$.}
\subsubsection{Splitting on subregions}
We first split the region $(0,\IF)^2$ into the following two parts: 
\BQNY
 \mcE_1:=[t_A-\theta_0, t_A+\theta_0]\times [s_A-\theta_0, s_A+\theta_0]\subset A. \  \ \
\mcE_2:=(0,\IF)^2\setminus \mcE_1,
\EQNY
where $\theta_0>0$ is some small constant which can be identified later on.
It follows from \eqref{eq:SelfSimi} that
\BQN \label{eq:Pu_bounds}
&&P_{0}(u):=\pk{\exists_{(t,s)\in   \mcE_1   }\    X_1(t)>\sqrt u (1+\mu_1 t),  X_2(s)>\sqrt u (1+\mu_2 s)} \nonumber \\
&&\le P(u)\le \pk{\exists_{(t,s)\in  \mcE_1   }\    X_1(t)>\sqrt u (1+\mu_1 t),  X_2(s)>\sqrt u (1+\mu_2 s)}  \\
&&\ \ \ \ \ \ \ \   \ \ \  \  + \pk{\exists_{(t,s)\in  \mcE_2 }\    X_1(t)>\sqrt u (1+\mu_1 t),  X_2(s)>\sqrt u (1+\mu_2 s)} \nonumber\\
&&=:P_{0}(u)+r_0(u)\nonumber
\EQN
Furthermore, we have, for all large $u$,
\BQN\label{eq:pr}
p(u)\le P_{0}(u)\le  p(u)+ r_1(u),
\EQN
where
\BQNY
&&p(u):=\pk{\exists_{(t,s)\in \Del^{(1)}_u\times \Del^{(2)}_u}\    X_1(t)>\sqrt u (1+\mu_1 t),  X_2(s)>\sqrt u (1+\mu_2 s)}, \\
&& r_1(u):=\pk{ \exists _{(t,s) \in \mcE_1\setminus \Del^{(1)}_u\times \Del_u^{(2)} } \    X_1(t)>\sqrt u (1+\mu_1 t),  X_2(s)>\sqrt u (1+\mu_2 s) },
\EQNY
with
\BQNY
&&\Del^{(1)}_u=\left[t_A-\frac{\ln(u)}{\sqrt{u}}, t_A+\frac{\ln(u)}{\sqrt{u}} \right],\ \Del^{(2)}_u=\left[s_A-\frac{\ln(u)}{\sqrt{u}}, s_A+\frac{\ln(u)}{\sqrt{u}} \right]. 
\EQNY
Next, we further split the rectangle $\Del^{(1)}_u\times \Del^{(2)}_u$ into smaller rectangles. To this end,
we denote, for any fixed $T, S>0$ 
\BQNY 
&&\Del^{(1)}_{j;u}=\Del^{(1)}_{j;u}(T)= [t_A+ j Tu^{-1}, t_A+( j+1) Tu^{-1}],\ \ -N^{(1)}_u  \le j \le N^{(1)}_u,\\
&&\Del^{(2)}_{l;u}=\Del^{(2)}_{l;u}(S)= [s_A+ l Su^{-1}, s_A+( l+1) Su^{-1}], \ \ -N^{(2)}_u  \le   l\le N^{(2)}_u,
\EQNY
where $N^{(1)}_u=\lfloor T^{-1} \ln(u) \sqrt{u}\rfloor$, $N^{(2)}_u=\lfloor S^{-1} \ln(u) \sqrt{u}\rfloor$
 (we denote by $\lfloor a \rfloor$ the smallest integer that is larger than $a$). Define
\begin{equation*} 
p_{j,l;u}=\pk{\exists_{ (t,s)\in\Del^{(1)}_{j;u}\times \Del_{l;u}^{(2)}} \ X_1(t)>
\sqrt{u}(1+ \mu_1 t ), X_2(s)>
\sqrt{u}(1+ \mu_2 s )}
\end{equation*}
and
\BQNY
&&p_{j,l_1,l_2;u}=\pk{\exists_{t\in \Del^{(1)}_{j;u}}   {X}_1(t)>\sqrt u (1+ \mu_1 t), \exists_{s\in \Del^{(2)}_{l_1;u}} {X}_2(s)>\sqrt u (1+ \mu_2 s), \exists_{s\in \Del^{(2)}_{l_2;u}}  {X}_2(s)>\sqrt u (1+ \mu_2 s)  }\\
&&\overline{p}_{j_1,j_2,l;u}=\pk{\exists_{t\in \Del^{(1)}_{j_1;u}}  {X}_1(t)>\sqrt u (1+ \mu_1 t), \exists_{t\in \Del^{(1)}_{j_2;u}}  {X}_1(t)>\sqrt u (1+ \mu_1 t), \exists_{s\in \Del^{(2)}_{l;u}} {X}_2(s)>\sqrt u (1+ \mu_2 s) }.
\EQNY
We have from the generalized Bonferroni's inequality (see Lemma \ref{Lem:Bonf} in Appendix \ref{a.B})
\BQN\label{eq:thetaT1}
p_1(u)\ge p(u) \ge p_2(u)-\Pi_1(u)-\Pi_2(u),
\EQN
where
\BQNY 
&&p_1(u)= \sum_{j=-N^{(1)}_u}^{N^{(1)}_u} \sum_{l=-N^{(2)}_u}^{N^{(2)}_u} p_{j,l;u},\ \
p_2(u) = \sum_{j=-N^{(1)}_u+1}^{N^{(1)}_u-1} \sum_{l=-N^{(2)}_u+1}^{N^{(2)}_u-1} p_{j,l;u},\\
&&\Pi_1(u) = \sum_{j=-N^{(1)}_u}^{N^{(1)}_u}\sum_{ -N^{(2)}_u\le l_1< l_2\le  N^{(2)}_u}p_{j,l_1,l_2;u}, \ \
\Pi_2(u) = \sum_{l=-N^{(2)}_u}^{N^{(2)}_u}\sum_{ -N^{(1)}_u\le j_1< j_2\le  N^{(1)}_u} \overline{ p}_{j_1,j_2,l;u}.
\EQNY

\subsubsection{Upper bounds and estimates }
In what follows, we shall derive upper bounds for   $r_0(u), r_1(u)$ in Lemma \ref{Lem:r01},  the exact asymptotics of $p_1(u), p_2(u)$ in Lemma \ref{Lem:PL} and asymptotic behaviour for $\Pi_1(u), \Pi_2(u)$ in Lemma \ref{Lem:Pi12}.
The \kkkk{proofs}  of the lemmas are displayed in Appendix \ref{a.B}. Recall that we assume
$-1< \rho<\hat\rho_1$.

\BEL \label{Lem:r01}
For any chosen small $\theta_0>0,$   we have, for all large $u$,
\BQNY 
&&r_0(u) \le    e^{- \frac{(\sqrt u -C_0)^2 }{2 } \widehat g }, \ \ \text{with}\ \ \widehat g =\inf_{(t,s)\in  \mcE_2 } g(t,s) > g_A(t_A, s_A),\\
&&r_1(u) \le  C_1  u^{3/2}   e^{- \frac{u }{2 } g_A(t_A, s_A) -  K_1(\ln(u))^2}
\EQNY
hold  for some constants  $C_0, C_1,K_1>0$ not depending on $u$.
\EEL
Below we discuss the asymptotics of $p_1(u), p_2(u)$. Define
\BQNY
\H(\mu;T)  :=  \int_{\R}e^{2\mu x_1} \pk{\exists_{t\in[0,T]}\  B_1(t)-\mu t> x_1}dx_1.
\EQNY

\BEL\label{Lem:PL} We have, as $u\to\infty$,
\begin{eqnarray*} 
  p_1(u) \sim   p_2(u) \sim
\frac{\H(\mu_1;T)\H(\mu_2-2\mu_1\rho;S)}{T S} \frac{1}{\mu_1 (\mu_2-2\mu_1\rho)}  e^{ - \frac{ g_A(t_A,s_A)  }{2} u   }  .
\end{eqnarray*}

\EEL

The last lemma is concerned with the asymptotic behaviour of $\Pi_{1}(u), \Pi_2(u)$. 

\BEL \label{Lem:Pi12}
It holds that
\BQNY
\limsup_{S\to\IF}\limsup_{T\to\IF}\lim_{u\to\IF}\frac{\Pi_{1}(u)}{\exp(-g_A(t_A, s_A)u/2)}=\limsup_{S\to\IF}\limsup_{T\to\IF}\lim_{u\to\IF}\frac{\Pi_{2}(u)}{\exp(-g_A(t_A, s_A)u/2)}=  0.
\EQNY
\EEL
\subsubsection{Asymptotics of $P(u)$}
By
Lemmas \ref{Lem:r01}, \ref{Lem:PL}, \ref{Lem:Pi12} applied to \eqref{eq:Pu_bounds} - \eqref{eq:thetaT1}
we obtain that
\[
P(u)
\sim
\lim_{S\to\IF}\lim_{T\to\IF} \frac{\H(\mu_1;T)\H(\mu_2-2\mu_1\rho;S)}{T S} \frac{1}{\mu_1 (\mu_2-2\mu_1\rho)}  e^{ - \frac{ g_A(t_A,s_A)  }{2} u   }
=
e^{ - \frac{ g_A(t_A,s_A)  }{2} u},
\]
where we used that, for any $\mu>0$
\BQN\label{eq:Hmu}
\H(\mu):=\lim_{T\to \IF}\frac{1}{T} \H(\mu;T)=\mu,
\EQN
see, e.g., \cite{DHJT18}.
Hence, using that
$g_A(t_A,s_A)=
4(\mu_2+(1-2\rho)\mu_1)$ (see (i)-(ii) of Proposition \ref{Lem:Optg})
we conclude the proof for scenario  $-1< \rho<\hat\rho_1$ in Theorem \ref{Thm:mu12}.
 \QED

\subsection{(ii) Scenario $\hat\rho_1<\rho<\hat\rho_2$}
\subsubsection{Splitting on subregions}

 \begin{figure}
 \vspace{0mm}
  \includegraphics[width=180mm, height=70mm]{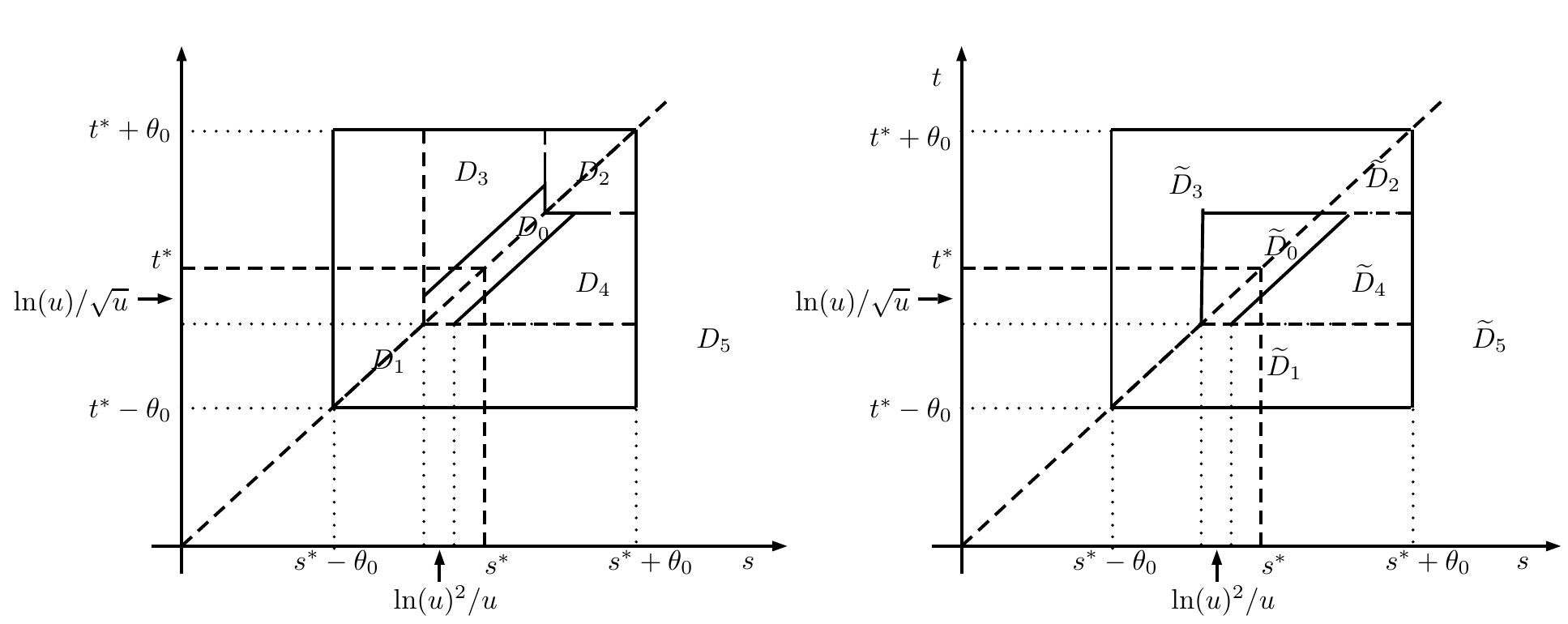}
  \caption{Partition of $(0,\IF)^2$: Left for $\hat\rho_1<\rho<\hat\rho_2$; right for $\rho=\hat\rho_1$}
\end{figure}

We split the region $(0,\IF)^2$ into five pieces as shown in Figure 1 (left). Namely, with some small $\theta_0>0$ and $u $ large,
let
\BQNY
&&D_0=\{(t,s): t^*- {\ln(u)}/{\sqrt{u}}\le t \le  t^*+ {\ln(u)}/{\sqrt{u}}, \ 0\le s-t \le \ln(u)^2/u \}\cup\\
&& \ \ \ \ \ \ \ \   \ \ \ \ \{(t,s): s^*- {\ln(u)}/{\sqrt{u}}\le s \le  s^*+ {\ln(u)}/{\sqrt{u}},\  0\le t-s\le \ln(u)^2/u \},\\
&&D_2= \{(t,s): t^*+ {\ln(u)}/{\sqrt{u}}\le t \le  t^*+ \theta_0, \ s^*+ {\ln(u)}/{\sqrt{u}}\le s \le  s^*+ \theta_0,  \},\\
&&D_3=\{(t,s): s^*- {\ln(u)}/{\sqrt{u}}\le s \le  s^*+ {\ln(u)}/{\sqrt{u}}, \ s+ {\ln(u)^2}/{u}\le t \le  t^*+ \theta_0,  \}, \\
&&D_4=\{(t,s): t^*- {\ln(u)}/{\sqrt{u}}\le t \le  t^*+ {\ln(u)}/{\sqrt{u}}, \ t+ {\ln(u)^2}/{u}\le s \le  s^*+ \theta_0,  \}, \\
&&D_1=[t^*-\theta_0, t^*+\theta_0]\times [s^*-\theta_0, s^*+\theta_0] \setminus (D_0 \cup D_2 \cup D_3 \cup D_4) , \\
&& D_5= (0,\IF)^2\setminus [t^*-\theta_0, t^*+\theta_0]\times [s^*-\theta_0, s^*+\theta_0] .
\EQNY
Clearly, we have  the following bounds
\BQN\label{eq:prrr}
p(u)\le P (u)\le  p(u)+ r_1(u)+r_2(u)+r_3(u),
\EQN
where
\BQNY
&&p(u):=\pk{\exists_{(t,s)\in D_0}\    X_1(t)>\sqrt u (1+\mu_1 t),  X_2(s)>\sqrt u (1+\mu_2 s)}, \\
&& r_1(u):=\pk{ \exists _{(t,s) \in D_5 } \    X_1(t)>\sqrt u (1+\mu_1 t),  X_2(s)>\sqrt u (1+\mu_2 s) },\\
&&r_2(u):=\pk{ \exists _{(t,s) \in D_1\cup D_2 } \    X_1(t)>\sqrt u (1+\mu_1 t),  X_2(s)>\sqrt u (1+\mu_2 s) },\\
&&r_3(u):=\pk{ \exists _{(t,s) \in D_3\cup D_4 } \    X_1(t)>\sqrt u (1+\mu_1 t),  X_2(s)>\sqrt u (1+\mu_2 s) }.
\EQNY
Next, we consider a further partition of $D_0$. Recall $\Del_{T,S}$ given in \eqref{eq:DelTS}. Denote, for any $T,S>0$ and $u>0,$
\BQNY 
&& \Del^{(1)}_{j;u}= \Del^{(1)}_{j;u}(T)= [t^*+ j Tu^{-1}, t^*+( j+1) Tu^{-1}],\ \ - N^{(1)}_u  \le j\le   N^{(1)}_u, \\
&&
  \Del^{(2)}_{l;u}=  \Del^{(2)}_{l;u}(S)= [ l Su^{-1}, ( l+1) Su^{-1}], \ \  1 \le  l\le   N_u^{(2)},
\EQNY
where $  N_u^{(1)}=\lfloor T^{-1} \ln(u) \sqrt{u}\rfloor$,  $  N_u^{(2)}=\lfloor S^{-1} \ln(u)^2\rfloor$.   Define further
\BQNY 
p_{j;u}&:=&\pk{\exists_{ ( t,s) \in (t^*+\frac{j T}{u},s^*+\frac{j T}{u})+u^{-1}\Del_{T,S}} \ X_1(t)>
\sqrt{u}(1+ \mu_1 t ), X_2(s)>
\sqrt{u}(1+ \mu_2 s )},\\
p_{j,l;u}&:=&\pk{\exists_{  t \in\Del^{(1)}_{j;u}, s-t\in \Del_{l;u}^{(2)}} \ X_1(t)>
\sqrt{u}(1+ \mu_1 t ), X_2(s)>
\sqrt{u}(1+ \mu_2 s )},\\
\overline{p}_{j,l;u}&:=&\pk{\exists_{  s \in\Del^{(1)}_{j;u}, t-s\in \Del_{l;u}^{(2)}} \ X_1(t)>
\sqrt{u}(1+ \mu_1 t ), X_2(s)>
\sqrt{u}(1+ \mu_2 s )},
\EQNY
and
\BQNY
&&q_{j_1,j_2;u}=\pk{
\begin{array}{cc}
\exists_{  ( t,s) \in (t^*+\frac{j_1 T}{u},s^*+\frac{j_1 T}{u})+u^{-1}\Del_{T,S}} \ X_1(t)>
\sqrt{u}(1+ \mu_1 t ), X_2(s)>
\sqrt{u}(1+ \mu_2 s ) \\
 \exists_{ ( t,s) \in (t^*+\frac{j_2 T}{u},s^*+\frac{j_2 T}{u})+u^{-1}\Del_{T,S}} \ X_1(t)>
\sqrt{u}(1+ \mu_1 t ), X_2(s)>
\sqrt{u}(1+ \mu_2 s )
 \end{array}
 }.
\EQNY
Thus, it follows from  the  Bonferroni's inequality that
\BQN\label{eq:thetaT}
\overline\Pi_1(u)+\Pi_1(u)+p_1(u)\ge p(u) \ge p_2(u)-\Pi_{21}(u)-\Pi_{22}(u),
\EQN
where
\BQNY 
&&p_1(u)= \sum_{j=-N_u^{(1)}}^{N_u^{(1)}}  p_{j;u},\ \ \ \
p_2(u) = \sum_{j=-N_u^{(1)}+1}^{N_u^{(1)}-1}  p_{j;u},\ \ \ \ \overline\Pi_1(u) = \sum_{j=-N_u^{(1)}}^{N_u^{(1)}}\sum_{ 1\le l \le  N_u^{(2)}} \overline p_{j,l;u},\\
&&\Pi_1(u) = \sum_{j=-N_u^{(1)}}^{N_u^{(1)}}\sum_{ 1\le l \le  N_u^{(2)}}p_{j,l;u}, \ \ \ \
\Pi_{21}(u):=\sum^{ N_u^{(1)}} _{ j_1=-N_u^{(1)}} \sum_{j_2>j_1+1} q_{j_1,j_2;u},\ \ \ \
\Pi_{22}(u):=\sum^{ N_u^{(1)}} _{ j_1=- N_u^{(1)}}   q_{j_1,j_1+1;u}.
\EQNY

\subsubsection{Upper bounds and estimates }
In what follows, we shall derive upper bounds for  $r_i(u),i =1,2,3$ in Lemma \ref{Lem:r123},
the exact asymptotics of $p_1(u), p_2(u)$ in Lemma \ref{Lem:PL2} and asymptotic behaviour for
$\overline \Pi_1(u), \Pi_1(u), \Pi_{21}(u),\Pi_{22}(u)$ in Lemma \ref{Lem:Pi112}.
The \kkkk{proofs}  of the lemmas are displayed in Appendix \ref{a.B}.


\BEL \label{Lem:r123}
For any chosen small $\theta_0>0,$   we have, for all large $u$,
\BQNY 
&&r_1(u) \le    e^{- \frac{(\sqrt u -C_1)^2 }{2 } \widehat g }, \ \ \ \widehat g =\inf_{(t,s)\in  D_5 } g(t,s) > g_L(t^*),\\
&&r_2(u) \le  C_2  u^{3/2}   e^{- \frac{u }{2 } g_L(t^*) -  K_2(\ln(u))^2},\\
&& r_3(u) \le  C_3  u^{3/2}   e^{- \frac{u }{2 } g_L(t^*) -  K_3(\ln(u))^2}
\EQNY
hold for some constants  $C_1, C_2, C_3, K_2, K_3>0$ not depending on $u$.

\EEL




\BEL\label{Lem:PL2} For any $T,S>0,$ we have, as $u\to\infty$,
\begin{eqnarray*} 
  p_1(u) \sim   p_2(u) \sim
\frac{\H(T,S)}{T}  \frac{\sqrt{t^*}}{2\sqrt{\pi(1-\rho)}} \  u^{-1/2}\ e^{ - \frac{u   g_L(t^*)  }{2}  }.
\end{eqnarray*}

\EEL

Below, we show, for any fixed $S>0$, the sub-additivity property of $\H(T,S)$ as a function of $T>0$.
\BEL\label{Lem:subadd}
Let $S>0$ be fixed, we have for any $T_1,T_2>0$
\BQNY
 \H(T_1+T_2,S)\le \H(T_1,S)+ \H(T_2, S)
\EQNY
and further, 
\BQNY
0< \frac{t^* \vk\mu^\top \Sigma_*^{-1} \vk b_*}{16 \prod_{i=1}^2(\Sigma_*^{-1} \vk b_*)_i}<\lim_{T\to\IF}\frac{1}{T}\H(T,S) = \inf_{T> 0}\frac{1}{T}\H(T,S) <\IF.
\EQNY

\EEL
The last lemma gives some asymptotic results for $\overline \Pi_1(u), \Pi_1(u), \Pi_{21}(u),\Pi_{22}(u)$.
\BEL\label{Lem:Pi112}
For any  $T>S>1$,
\BQNY
\lim_{u\to\IF}\frac{\max (\Pi_{1}(u),\overline\Pi_1(u))}{u^{-1/2}\exp(-g_L(t^*)u/2)}  \le C_0  \lfloor S \rfloor   \sum_{ l\ge 1 } e^{-K_0 lS},
\EQNY
\BQNY
\lim_{u\to\IF}\frac{\Pi_{21}(u)}{u^{-1/2}\exp(-g_L(t^*)u/2)}  \le C_1(S)  \lfloor T \rfloor   \sum_{ l\ge 1 } e^{-K_1 l T},
\EQNY
and
\BQNY
\lim_{u\to\IF}\frac{\Pi_{22}(u)}{u^{-1/2}\exp(-g_L(t^*)u/2)}  =\frac{\sqrt{t^*}}{2\sqrt{\pi(1-\rho)}} \ \LT(\frac{2\H(T,S)}{T}- \frac{\H(2T,S)}{T} \RT) ,
\EQNY
where $C_0, K_0, K_1>0$ are three constants which do not dependent on $T,S,u$, and $C_1(S)$ does not dependent on $T,u$.

\EEL

\subsubsection{Asymptotics of $P(u)$}
Combining \eqref{eq:prrr}-\eqref{eq:thetaT} and the  results in Lemmas \ref{Lem:r123}, \ref{Lem:PL2} and \ref{Lem:Pi112}, yields that, for any large $T_1,T_2, S_1,S_2$ such that $S_i< T_i, i=1,2$,
\BQNY
&& \frac{\sqrt{t^*}}{2\sqrt{\pi(1-\rho)}} \frac{\H(T_1,S_1)}{T_1}+ 2 C_0  \lfloor S_1 \rfloor   \sum_{ l\ge 1 } e^{-K_0 lS_1} \\
&&\ge \limsup_{u\to\IF}\frac{P(u)}{u^{-1/2}\exp(-g_L(t^*)u/2)} \ge \liminf_{u\to\IF}\frac{P(u)}{u^{-1/2}\exp(-g_L(t^*)u/2)}\\
&&\ge \frac{\sqrt{t^*}}{2\sqrt{\pi(1-\rho)}} \frac{\H(T_2,S_2)}{T_2} -
C_1(S_2)\lfloor T_2  \rfloor \
 \sum_{k\ge 1}  e^{-K_1 kT_2} - \frac{\sqrt{t^*}}{2\sqrt{\pi(1-\rho)}}  \LT(\frac{2\H(T_2,S_2)}{T_2}-\frac{\H(2T_2,S_2)}{T_2}\RT).
\EQNY
Letting first $T_2\to\IF$ and then $S_2\to\IF$, we have from the above formula, \eqref{eq:Hmu} and  \nelem{Lem:subadd} that
\BQNY
\lim_{S\to\IF}\lim_{T\to\IF}  \frac{\H(T,S)}{T}\in (0,\IF).
\EQNY
The proof for scenario $\hat\rho_1<\rho<\hat\rho_2$ in Theorem \ref{Thm:mu12} follows then by letting
$T_1\to\IF$ and then $S_1\to\IF$,
and (iv) of Proposition \ref{Lem:Optg}.
 \QED

\subsection{(iii) Scenario $\rho=\hat\rho_1$}
Since the idea of the proof of this case is {similar to  that} of scenarios (i) and (ii),
we present only main steps.
We split the region $(0,\IF)^2$ into five pieces as shown in Figure 1 (right).
Namely, with some small $\theta_0>0$ and $u $ large, let
\BQNY
&&\widetilde D_0=\{(t,s): t^*- {\ln(u)}/{\sqrt{u}}\le t \le  t^*+ {\ln(u)}/{\sqrt{u}}, \ 0\le s-t \le \ln(u)^2/u \}\cup\\
&& \ \ \ \ \ \ \ \   \ \ \ \ \{(t,s): s^*- {\ln(u)}/{\sqrt{u}}\le s \le  s^*+ {\ln(u)}/{\sqrt{u}},\  s< t \le t^*+{\ln(u)}/{\sqrt{u}} \}=:\widetilde D_{0B} \cup \widetilde D_{0A},\\
&&\widetilde D_3=\{(t,s): s^*- \theta_0 \le s \le  s^*+ \theta_0,\  s<t \le t^*+\theta_0\}\setminus \widetilde D_{0A},\\
&&\widetilde D_1=D_1\cap \overline B,\ \ \ \widetilde D_2=D_2\cap \overline B,\ \ \  \widetilde D_4=D_4,\ \ \ \widetilde D_5=D_5.
\EQNY
Clearly, we have  the following bounds
\BQN\label{eq:ptr}
p(u)\le P (u)\le  p(u)+ \tilde r_0(u)+\tilde r_1(u)+\tilde r_2(u),
\EQN
where
\BQNY
&&p(u):=\pk{\exists_{(t,s)\in \widetilde D_{0A}}\    X_1(t)>\sqrt u (1+\mu_1 t),  X_2(s)>\sqrt u (1+\mu_2 s)}, \\
&& \tilde r_0(u):=\pk{ \exists _{(t,s) \in  \widetilde D_{0B}} \    X_1(t)>\sqrt u (1+\mu_1 t),  X_2(s)>\sqrt u (1+\mu_2 s) },\\
&& \tilde r_1(u):=\pk{ \exists _{(t,s) \in  \widetilde D_1\cup \widetilde D_2\cup \widetilde D_4\cup\widetilde D_5} \    X_1(t)>\sqrt u (1+\mu_1 t),  X_2(s)>\sqrt u (1+\mu_2 s) },\\
&&\tilde r_2(u):=\pk{ \exists _{(t,s) \in \widetilde D_3 } \    X_1(t)>\sqrt u (1+\mu_1 t),  X_2(s)>\sqrt u (1+\mu_2 s) }.
\EQNY
Similar arguments as used in scenarios (i), (ii) give that
\BQN \label{eq:puxx}
\limit{u} \frac{p(u)}{ \exp\LT(-  g_A(t_A, s_A) u/  2     \RT) } =  \frac{\mu_1^{3/2} (\mu_2-2\mu_1\rho)^2 }{2\pi \sqrt{\mu_2-2(\mu_1+\mu_2) \rho +3\mu_1 \rho ^2} }
\int_{\R}\int_{x_2}^\IF e^{\frac{- (a_1x_1^2-2a_2 x_1x_2+a_3 x_2^2)}{4}} \,dx_1 dx_2,
\EQN
and
\BQN \label{eq:r0t}
 \limit{u} \frac{\tilde r_0(u) }{  u^{-1/2}  \exp\LT(-  g_A(t_A, s_A) u/  2     \RT) } \le  \frac{\widetilde\H \sqrt{t^*}}{2\sqrt{\pi(1-\rho)}},
\EQN
and the asymptotically negligibility of $\tilde r_1(u), \tilde r_2(u)$. Note that in proving the bound for $\tilde r_2(u)$, in addition to \eqref{eq:gt0s0} as in the proof of Lemma \ref{Lem:r01}, we also need the fact that (for $t>s$)
\BQNY
g_A(t_A+t,s_A+s)&\ge& g_A(t_A,s_A)+\frac{a_1}{2}(1-\vn)\Bigg( \LT( t+\frac{\mu_2-2\mu_1\rho}{\mu_1}\rho s\RT)^2\nonumber\\
&&\ \ \ \ \ \ \ \ \ \ \ \ \ \ \ \ \ \ \ \ \ \ \ \ \ \ \ \ \ \ \  \ \ +\LT(\frac{(\mu_2-2\mu_1\rho)^3(\mu_2-2(\mu_1+\mu_2) \rho +3\mu_1 \rho ^2)}{\mu_1^3}\RT)  s^2 \Bigg).
\EQNY
Consequently, the claim follows by formulas \eqref{eq:ptr}-\eqref{eq:r0t} and the asymptotically negligibility
of $\tilde r_1(u), \tilde r_2(u)$. This completes the proof of scenario $\rho=\hat\rho_1$ in Theorem \ref{Thm:mu12}.
 \QED
\subsection{(iv) Scenario $\hat\rho_2\le\rho<1$}
\label{sc.4}
First note that
\[e^{-2\mu_2 u}=\pk{\sup_{s\ge 0} (X_2(s)-\mu_2 s)>u} \ge  P(u)\\
\ge  \pk{\exists_{t\ge  0} \   X_1(t)-\mu_1 t>u, \   X_2(t)-\mu_2 t>u}=:\pi(u).
\]
Furthermore, the exact asymptotics for $\pi(u)$  has been discussed in  Corollary 4.3 in \cite{Ji18} (where we  take $r=0$).  Thus, we have, for $\rho=\hat\rho_2$,
\BQNY
\pi(u)\sim  
\frac{1}{2} e^{-2\mu_2 u}, \ \ \ u\to\IF,
\EQNY
 and for  $\hat\rho_2<\rho<1$,
 \BQNY
\pi(u)\sim   e^{-2\mu_2 u},\ \ \ u\to\IF.
\EQNY
Therefore, the claims in scenario $\hat\rho_2<\rho<1$ of Theorem \ref{Thm:mu12}
and $\hat\rho_2=\rho$ in (a) of Remark \ref{Rem:main1} follow.
\QED



\section{Proof of Theorem \ref{Thm:mu}} \label{proofmain2}

For $\mu_1 = \mu_2=\mu$, we have that $\hat\rho_1=0, \hat\rho_2=1$.
The case $\rho=0$ follows from (\ref{eq:Prho0}).
Thus the interesting scenarios include
(i) $-1< \rho< 0$ and (ii) $0<\rho<1$. The claim for (ii) $0<\rho<1$ follows directly from (iii) in Theorem \ref{Thm:mu12}, with $t^*=1/\mu$. Next, we shall focus on the proof for (i) $-1< \rho< 0$. The proof goes with the same arguments as in the proof of scenario (i) in Theorem \ref{Thm:mu12}, but now there are two minimizers of the function $g(t,s), (t,s)\in(0,\IF)^2,$ namely, $(t_A,s_A)=(t_0,s_0)\in A, (t_B,s_B)=(s_0,t_0)\in B$, with
$
t_0=   \frac{ 1-2\rho}{  \mu},  s_0 =  \frac{1}{(1-2\rho)\mu}.
$

We first split the region $(0,\IF)^2$ into three parts. Namely, with some small $\theta_0>0$, let
\BQNY
&& \mcE_{11}=[t_0-\theta_0, t_0+\theta_0]\times [s_0-\theta_0, s_0+\theta_0]\subset A, \\
&& \mcE_{12}=[s_0-\theta_0, s_0+\theta_0]\times [t_0-\theta_0, t_0+\theta_0] \subset B,\\
&& \mcE_{2}=(0,\IF)^2\setminus (\mcE_{11}\cup \mcE_{12}).
\EQNY
\COM{ 
It follows that
\BQN \label{eq:Pu_bounds}
&&\pk{\exists_{(t,s)\in  (D_1\cup D_2)  }\    X_1(t)>\sqrt u (1+\mu_1 t),  X_2(s)>\sqrt u (1+\mu_2 s)} \nonumber \\
&&\le P(u)\le \pk{\exists_{(t,s)\in  (D_1\cup D_2)  }\    X_1(t)>\sqrt u (1+\mu_1 t),  X_2(s)>\sqrt u (1+\mu_2 s)}  \\
&&\ \ \ \ \ \ \ \   \ \ \  \  + \pk{\exists_{(t,s)\in  D_3  }\    X_1(t)>\sqrt u (1+\mu_1 t),  X_2(s)>\sqrt u (1+\mu_2 s)} \nonumber
\EQN
}
As in  the proof of scenario (i) of Theorem \ref{Thm:mu12}, the main contribution of the asymptotics comes from $\mcE_{11}\cup \mcE_{12}$.
Note further that
\BQNY
\lefteqn{\pk{\exists_{(t,s)\in  (\mcE_{11}\cup \mcE_{12})  }\    X_1(t)>\sqrt u (1+\mu_1 t),  X_2(s)>\sqrt u (1+\mu_2 s)}} \nonumber\\
& =& \pk{\exists_{(t,s)\in   \mcE_{11}  }\    X_1(t)>\sqrt u (1+\mu_1 t),  X_2(s)>\sqrt u (1+\mu_2 s)}  \\
&&+ \pk{\exists_{(t,s)\in   \mcE_{12}}\    X_1(t)>\sqrt u (1+\mu_1 t),  X_2(s)>\sqrt u (1+\mu_2 s)}  \nonumber\\
&&-\pk{
\begin{array}{cc}
\exists_{(t_1,s_1)\in   \mcE_{11} } \  X_1(t_1)>\sqrt u (1+\mu_1 t_1),  X_2(s_1)>\sqrt u (1+\mu_2 s_1)\\
 \exists_{(t_2,s_2)\in  \mcE_{12}}    \    X_1(t_2)>\sqrt u (1+\mu_1 t_2),  X_2(s_2)>\sqrt u (1+\mu_2 s_2)
\end{array} }\nonumber\\
&=:&P_{\theta_0,1}(u)+P_{\theta_0,2}(u)-P_{\theta_0,0}(u).
\EQNY
By symmetric property of the model we know that $P_{\theta_0,1}(u)=P_{\theta_0,2}(u)$. Next, we show in Lemma \ref{Lem:Pr2} that $P_{\theta_0,0}(u)$ is asymptotically negligible compared with $P_{\theta_0,1}(u)$. 
The proof of it is displayed in Appendix \ref{a.B}.

\BEL \label{Lem:Pr2}
For any chosen small  $\theta_0>0,$ we have for all large $u$
\BQNY 
P_{\theta_0,0}(u) \le    e^{-\frac{( \sqrt u   -C_0)^2g_A(t_A,s_A)}{2\sigma_0^2}}
\EQNY
holds  for some constant  $ C_0>0, \sigma_0^2\in(0,1)$ which do not depend on $u$.
\EEL
The rest of the proof is the same as those in  the proof of scenario (i) in Theorem \ref{Thm:mu12}, and thus omitted. This completes the proof. \QED

\begin{appendix}


\section{Proofs of Lemmas \ref{Lem:r01}-\ref{Lem:Pr2} }\label{a.B}

In this section we give proofs of Lemmas \ref{Lem:r01}-\ref{Lem:Pr2} that are
the building blocks of the proofs of Theorems \ref{Thm:mu12} and \ref{Thm:mu}.

We begin with the analysis of the local behaviour of  function $g(t,s), (t,s)\in(0,\IF)$
at its minimizer in scenarios \jj{(i)--(iv)} of Proposition \ref{Lem:Optg}, respectivelly.   

\BEL\label{Lem:Taylor} Assume that $\mu_1<\mu_2$. We have
\begin{itemize}
\item[(i).] If  $-1< \rho  <\hat \rho_1$, then as $(t,s)\to(0,0)$,
\BQNY 
g(t_A+t,s_A+s)=g_A(t_A,s_A)+\frac{a_1}{2}t^2\oo - a_2 ts\oo +\frac{a_3}{2} s^2\oo,
\EQNY
where, with $h(\rho):=\mu_2-2(\mu_1+\mu_2)\rho+3\mu_1\rho^2>0,$
\BQNY
 a_1:=\frac{2\mu_1^3(\mu_2-2\mu_1\rho)}{h(\rho)}>0, \ \
 a_2:=\frac{-2\rho\mu_1^2(\mu_2-2\mu_1\rho)^2}{h(\rho)},  \ \
 a_3: =\frac{2(\mu_2-2\mu_1\rho)^4(1-2\rho)}{h(\rho)}>0.
\EQNY

\item[(ii).] If $ \hat \rho_1  < \rho< \hat \rho_2 $, then
\begin{itemize}
\item{(ii.1),} as $(t,s)\to(0,0)$, with $s<t$ (i.e., $(t^*+t,s^*+s)\in A$),
\BQNY 
g(t^*+t,s^*+s)=g_L(t^*)+b_1( t-s )  \oo + \frac{ c_1 }{2}s^2\oo,
\EQNY
where
\BQNY
b_1&:=&\frac{( \rho-1-2\rho^2)+2\rho(\mu_2-\mu_1\rho)s^* +(1+\rho)\mu_1^2{s^*}^2}{(1-\rho)(1+\rho)^2{s^*}^2}>0,\\
c_1&:=& \frac{2 }{ {s^*}^3}\LT(1+\frac{\rho^2(\rho(1-\rho)-(\mu_2-\mu_1\rho) s^*)^2}{(1-\rho^2)^3}\RT)>0.
\EQNY
\item{(ii.2),} as $(t,s)\to(0,0)$, with $s>t$ (i.e., $(t^*+t,s^*+s)\in B$),
\BQNY 
g(t^*+t,s^*+s)=g_L(t^*)+b_2( s-t )  \oo + \frac{ c_2 }{2}t^2\oo,
\EQNY
where
\BQNY
b_2&:=&\frac{( \rho-1-2\rho^2)+2\rho(\mu_1-\mu_2\rho)t^* +(1+\rho)\mu_2^2{t^*}^2}{(1-\rho)(1+\rho)^2{t^*}^2}>0,  \\
c_2&:=& \frac{2 }{ {t^*}^3}\LT(1+\frac{\rho^2(\rho(1-\rho)-(\mu_1-\mu_2\rho) t^*)^2}{(1-\rho^2)^3}\RT)>0.
\EQNY
\item{(ii.3),} as $(t,s)\to(0,0)$, with $s=t$ (i.e., $(t^*+t,s^*+s)\in L$),
\BQNY 
g(t^*+t, s^*+t)=g_L(t^*)+\frac{b_0}{2}t^2\oo,
\EQNY
where
$b_0:=\frac{4}{(1+\rho){t^*}^3}.
$
\end{itemize}

\item[(iii).]  If $\rho = \hat \rho_1$ (in this case $t_A=s_A=t^*=s^*$), then
  \begin{itemize}
\item{(iii.1),} as $(t,s)\to(0,0)$, with $s<t$,
\BQNY
g(t_A+t,s_A+s)=g_A(t_A,s_A)+\frac{a_1}{2}t^2\oo - a_2 ts\oo +\frac{a_3}{2} s^2\oo,
\EQNY
\item{(iii.2),} as $(t,s)\to(0,0)$, with $s>t$,
\BQNY 
g(t^*+t,s^*+s)=g_L(t^*)+b_2( s-t )  \oo + \frac{ c_2 }{2}t^2\oo.
\EQNY
\item{(iii.3),} as $(t,s)\to(0,0)$, with $s=t$,
\BQNY 
g(t^*+t, s^*+t)=g_L(t^*)+\frac{b_0}{2}t^2\oo.
\EQNY

\end{itemize}


\end{itemize}

\EEL

\LLc{The proof of Lemma \ref{Lem:Taylor} is tedious but  only involves basic calculations using Taylor expansion, and thus it is omitted.
}


Next we present below a generalized version of the Bonferroni's inequality. The proof can be found in, e.g., \cite{HJ14c}.

\BEL\label{Lem:Bonf}
  Let $(\Omega, \mathcal{F}, \mathbb{P})$ be a probability space and $A_1,\cdots,A_n$ and $B_1,\cdots,B_m$ be $n+m$ events in $\mathcal{F}$ with $n,m\ge 2$. Then
\BQNY 
&&\sum_{k=1}^n\sum_{l=1}^m\pk{A_k\cap B_l}\ge \pk{\underset{l=1,\ldots, m}{\underset{k=1,\ldots,n} \bigcup }(A_k\cap B_l)}\ge\sum_{k=1}^n\sum_{l=1}^m\pk{A_k\cap B_l}\nonumber\\
&&\ \ \ \qquad -\sum_{k=1}^n\sum_{1\le l_1<l_2\le m}\pk{A_k\cap B_{l_1}\cap B_{l_2}}
-\sum_{l=1}^m\sum_{1\le k_1<k_2\le n}\pk{ A_{k_1}\cap A_{k_2}\cap B_l}.
\EQNY
 \EEL

 \subsection{Proof of Lemma \ref{Lem:r01}}
  Let $T_0>0$ be a  fixed large  constant (will be determined later). It is easily seen that
\BQNY
r_0(u)&\le& \pk{ \exists _{(t,s) \in [0,T_0]^2\setminus \mcE_1}  \    X_1(t)>\sqrt u (1+\mu_1 t),  X_2(s)>\sqrt u (1+\mu_2 s) }\\
&&+\pk{ \exists _{t\ge T_0}   \    X_1(t)>\sqrt u (1+\mu_1 t) }+\pk{\exists _{s\ge T_0}   \  X_2(s)>\sqrt u (1+\mu_2 s) }.
\EQNY
Next we consider upper bounds for each term on the right-hand side.
\LLc{According to \jj{Lemma 5} of \cite{DJT19a}, 
for any fixed $t,s$, there exists a unique index set
$$
I(t,s) \subseteq \{1,2\}
$$
such that
\BQN \label{eq:gstI}
g(t,s)=(1+\mu_1 t, 1+\mu_2 s)_{I(t,s)} \ (\Sigma_{ts})_{{I(t,s)},{I(t,s)}}^{-1} \  (1+\mu_1 t, 1+\mu_2 s)_{I(t,s)}^\top,
\EQN
and
\BQN\label{eq:SigII}
(\Sigma_{ts})_{{I(t,s)},{I(t,s)}}^{-1} \  (1+\mu_1 t, 1+\mu_2 s)_{I(t,s)}^\top > \vk 0_{I(t,s)}.
\EQN
Thus,}
\BQN \label{eq:DD}
&& \pk{ \exists _{(t,s) \in [0,T_0]^2\setminus  \mcE_1}  \    X_1(t)>\sqrt u (1+\mu_1 t),  X_2(s)>\sqrt u (1+\mu_2 s) }\nonumber\\
&& \le \pk{ \exists _{(t,s) \in [0,T_0]^2\setminus  \mcE_1}  \    (1+\mu_1 t, 1+\mu_2 s)_{I(t,s)}\ (\Sigma_{ts})_{{I(t,s)},{I(t,s)}}^{-1} \ ( X_1(t),X_2(s))_{I(t,s)}^\top >\sqrt u g(t,s)}\\
&&= \pk{ \exists _{(t,s) \in [0,T_0]^2\setminus  \mcE_1}  \   \frac{ Z(t,s)}{g(t,s)}>\sqrt u },\nonumber
\EQN
where
\BQN\label{eq:Z}
Z(t,s):= (1+\mu_1 t, 1+\mu_2 s)_{I(t,s)}\ (\Sigma_{ts})_{{I(t,s)},{I(t,s)}}^{-1} \ ( X_1(t),X_2(s))_{I(t,s)}^\top .
\EQN
Note that 
\BQN\label{eq:Var_Zg}
\mathrm{Var}\LT(\frac{Z(t,s)}{g(t,s)}\RT)=\frac{1}{g(t,s)}.
\EQN
\LLc{
In order to apply the Borell-TIS inequality, we first show that 
$$
\limsup_{(t,s)\to(t^{(b)},s^{(b)})}\frac{\abs{Z(t,s)}}{g(t,s)}<\IF,\ \ \ \ \text{almost\ surely}
$$
holds for any $(t^{(b)},s^{(b)})$ on the boundary $\{(t,s): t\ge0, s=0\}\cup\{(t,s): t=0, s\ge 0\}$.\\
In fact, if the above does not hold for some boundary point $(t^{(b)},s^{(b)})$, then  for any $M>0$ there   exist a sequence $\{(t_k,s_k)\}_{k=1}^\IF$ and some measurable set $E$ such that $(t_k,s_k) \to (t^{(b)},s^{(b)})$, $\pk{E}>0$ and
$$
\frac{\abs{Z(t_k,s_k)}}{g(t_k,s_k)}\ge M \ \ \ \text{on}\ E
$$
for all large enough $k$. Then we have
\BQN\label{eq:Var_LB}
\mathrm{Var}\LT(\frac{Z(t_k,s_k)}{g(t_k,s_k)} \RT)\ge M^2 \pk{E}>0.
\EQN
On the other hand, by \jj{Lemma 6} of \cite{DJT19a} we have $g(t,s)=g_3(t,s)$ for all $(t,s) \in  \{(t,s): t\ge0, s=0\}\cup\{(t,s): t=0, s\ge 0\}$, and thus by \eqref{eq:Var_Zg} and \eqref{eq:g32} we have
$\lim_{k\to\IF} \mathrm{Var}\LT(\frac{Z(t_k,s_k)}{g(t_k,s_k)} \RT) =0$. This is a contradiction with \eqref{eq:Var_LB}.  Therefore,  $\frac{ Z(t,s)}{g(t,s)}, (t,s) \in [0,T_0]^2\setminus  \mcE_1 $ is almost surely bounded.
}
Consequently, by the Borell-TIS inequality (see, e.g., \cite{AdlerTaylor}) we have, for any   fixed small constant  $\theta_0>0$
\BQNY 
\pk{ \exists _{(t,s) \in [0,T_0]^2\setminus   \mcE_1 }  \  \frac{Z(t,s)}{g(t,s)}>\sqrt u }\le e^{-\frac{(\sqrt u -C_0)^2}{2} \widehat g}
\EQNY
holds for all $u$ such that
$$\sqrt u>C_0:=\E{\sup_{(t,s)\in [0,T_0]^2\setminus   \mcE_1  } \frac{ Z(t,s)}{g(t,s)}}.$$

Moreover, since $X_i$ is the standard Brownian motion,
\BQNY
\lim_{t\to\IF}\frac{X_i(t)}{1+\mu_i t} = 0\ \  \ \ \mathrm{almost\  surely},
\EQNY
showing that the random process $\frac{X_i(t)}{1+\mu_i t}, t\ge T_0$ has almost surely bounded sample paths on $[T_0,\IF)$. Again by the Borell-TIS inequality
\BQNY 
\pk{ \exists _{t\ge T_0}   \    X_i(t)>\sqrt u (1+\mu_i t) }\le e^{-\frac{(\sqrt u -C_i)^2}{2}\frac{(1+\mu_i T_0)^2}{T_0}}
\EQNY
holds for all $\sqrt u>C_i:=\E{\sup_{t\in[T_0,\IF)} \frac{X_1(t)}{1+\mu_i t}}$.  Since for all large enough $T_0$  it holds that
$
\frac{(1+\mu_i T_0)^2}{T_0} >\widehat g,
$
the claim for $r_0(u)$ is established.

Below we consider $r_1(u)$. Since $(t_A,s_A)\in A$,  we have from Proposition \ref{Lem:Optg} that for any chosen small $\theta_0$
\BQNY
g(t,s)=g_A(t,s), \ \ (t,s)\in\mcE_1 \subset A,
\EQNY
and further   (cf. \eqref{eq:Z})
\BQNY
Z(t,s)= (1+\mu_1 t, 1+\mu_2 s)\ \Sigma_{ts}^{-1} \ ( X_1(t),X_2(s))^\top  =: h_1(t,s) X_1(t) + h_2(t,s) X_2(s), \ \ (t,s)\in\mcE_1,
\EQNY
with
\BQNY
h_1(t,s) =\frac{(1+\mu_1 t)s-\rho  s  (1+\mu_2 s))}{ts-\rho^2    s ^2},\ \ \
h_2(t,s) =\frac{(1+\mu_2 s)t-\rho     s (1+\mu_1 t))}{ts-\rho^2   s ^2}.
\EQNY
Thus, similarly to \eqref{eq:DD} we conclude that
\BQN \label{eq:Upp22}
r_1(u)\le \pk{ \exists _{(t,s) \in\mcE_1 \setminus \Del^{(1)}_u \times \Del_u^{(2)}  }  \   \frac{ Z(t,s)}{g_A(t,s)}>\sqrt u }.
\EQN
Since $
h_1(t,s), h_2(t,s), g_A(t,s),  (t,s)\in \mcE_1$
are all smooth functions and
\BQNY
\E{(X_i(t_1)-X_i(t_2))^2}=\abs{t_1-t_2}, \ \ i=1,2
\EQNY
one can check that, for all $(t_1,s_1), (t_2,s_2)\in \mcE_1,$
\BQNY
\E{\LT(\frac{Z(t_1,s_1)}{g_A(t_1,s_1)}-\frac{Z(t_2,s_2)}{g_A(t_2,s_2)}\RT)^2}\le \mathrm{Const}\cdot (\abs{t_1-t_2}+\abs{s_1-s_2}).
\EQNY
Therefore, an application of the Piterbarg's inequality  in \cite{KEP2015}[Lemma 5.1]
(see also \cite{Pit96}[Theorem 8.1] or \cite{MR3493177}[Theorem 3]) yields that
\BQN \label{eq:Upp22}
r_1(u)\le \pk{ \exists _{(t,s) \in \mcE_1 \setminus \Del^{(1)}_u \times \Del_u^{(2)}  }  \   \frac{ Z(t,s)}{g_A(t,s)}>\sqrt u }\le C_3 u^{3/2} e^{-\frac{  u }{2} \widetilde g_u},
\EQN
where $C_3>0$ is some constant which does not depend on $u$ and
\BQNY
\widetilde g_u : =  \inf_{(t,s)\in \mcE_1 \setminus \Del^{(1)}_u \times \Del_u^{(2)} } g_A(t,s).
\EQNY
Moreover, we have from (i) of Lemma \ref{Lem:Taylor} that for all $(t_A+t,s_A+s)\in \mcE_1$
\BQN\label{eq:gt0s0}
g_A(t_A+t,s_A+s)&\ge& g_A(t_A,s_A)+\frac{a_1}{2}(1-\vn)\Bigg( (1-\rho^2) t^2+\LT(\rho t+\frac{\mu_2-2\mu_1\rho}{\mu_1}s\RT)^2\nonumber\\
&&\ \ \ \ \ \ \ \ \ \ \ \ \ \ \ \ \ \ \ \ \ \ \ \ \ \ \ \ \ \ \  \ \ +\LT(\frac{\mu_2-2\mu_1\rho}{\mu_1}\RT)^2\LT(\frac{\mu_2-2\mu_1\rho}{\mu_1} (1-2\rho)-1\RT) s^2 \Bigg)
\EQN
holds with some small $\vn>0$, where for all $-1< \rho<\hat\rho_1$
{(see \kkk{also} the proof of (b).(i) in Lemma 9 of \cite{DJT19a} for $\rho>0$)}
\BQNY
\frac{\mu_2-2\mu_1\rho}{\mu_1} (1-2\rho)-1>0.
\EQNY
 Thus
\BQNY
\widetilde g_u \ge
  g_A(t_A,s_A) +\frac{a_1}{2}(1-\vn)\min\LT((1-\rho^2), \LT(\frac{\mu_2-2\mu_1\rho}{\mu_1}\RT)^2\LT(\frac{\mu_2-2\mu_1\rho}{\mu_1} (1-2\rho)-1\RT) \RT)\frac{(\ln(u))^2}{u}.
\EQNY

Inserting the above to \eqref{eq:Upp22} completes the proof.
\QED

\subsection{Proof of Lemma \ref{Lem:PL}}

\def \cu {\vk c_{\vk j;u}}
We first analyze the  summand $p_{j,l;u}$. We set
\BQN\label{eq:bba}
\vk b_{j,l;u}=(a_{j;u}, b_{l;u})^\top, \ \ \ \ a_{j;u}=1+\mu_1 (t_A+\frac{jT}{u}), \ \ b_{l;u}=1+\mu_2 (s_A+\frac{lS}{u}).
\EQN
It follows that
\begin{eqnarray}\label{eq:pjl_1}
p_{j,l;u}&=&\pk{ \underset{s\in[0,S]}{\exists_{t \in [0, T]}}
\begin{array}{ccc}
 X_1(t_A+\frac{jT}{u}+\frac{t}{u})> a_{j;u}
\sqrt{u} +\frac{\mu_1}{\sqrt u} t \\
 X_2(s_A+\frac{lS}{u}+\frac{s}{u})> b_{l;u}
\sqrt{u} +\frac{\mu_2}{\sqrt u} s
\end{array}
} \nonumber\\
&=&
\pk{ \underset{s\in[0,S]}{\exists_{t \in [0, T]}}
\begin{array}{ccc}
X_1(t_A+\frac{jT}{u})+ X_1(t_A+\frac{jT}{u}+\frac{t}{u})- X_1(t_A+\frac{jT}{u})> a_{j;u}
\sqrt{u} +\frac{\mu_1}{\sqrt u} t \\
 X_2(s_A+\frac{lS}{u})+ X_2(s_A+\frac{lS}{u}+\frac{s}{u})- X_2(s_A+\frac{lS}{u})> b_{l;u}
\sqrt{u} +\frac{\mu_2}{\sqrt u} s
\end{array}
} .
\end{eqnarray}
Since $(t_A+\frac{jT}{u},  s_A+\frac{lS}{u})\in A$ for all large $u,$ the covariance matrix of $\vk Z_{j,l;u}:=(X_1(t_A+\frac{jT}{u}),  X_2(s_A+\frac{lS}{u}))^\top$ is given by
$$
\  \Sigma_{j,l;u}=
\left(
   \begin{array}{cc}
     t_A+\frac{jT}{u}  & \rho\  (s_A+\frac{lS}{u})\\
     \rho\  (s_A+\frac{lS}{u}) & s_A+\frac{lS}{u} \\
   \end{array}
 \right).
$$
Thus, the density function of $\vk Z_{j,l;u}$ is given by
\BQNY
\phi_{\SI_{j,l;u}}( \vk w)=\frac{1}{\sqrt{(2\pi)^{2}\abs{\SI_{j,l;u}}}}\exp\LT(-\frac{1}{2}\vk w^\top(\Sigma_{j,l;u})^{-1}  \vk w  \RT),\ \ \ \vk w=(w_1,w_2)^\top.
\EQNY
By conditioning on the value of $\vk Z_{j,l;u}$ we rewrite \eqref{eq:pjl_1} as
\BQNY
p_{j,l;u}=\int_{\R^{2}}\phi_{\SI_{ j,l;u} }( \vk{w})
\pk{ \underset{s\in[0,S]}{\exists_{t \in [0, T]}}
\begin{array}{ccc}
X_1(t_A+\frac{jT}{u}+\frac{t}{u})- X_1(t_A+\frac{jT}{u})> a_{j;u}
\sqrt{u} +\frac{\mu_1}{\sqrt u} t -w_1\\
X_2(s_A+\frac{lS}{u}+\frac{s}{u})- X_2(s_A+\frac{lS}{u})> b_{l;u}
\sqrt{u} +\frac{\mu_2}{\sqrt u} s -w_2
\end{array}
\Bigg{|} \vk Z_{j,l;u}=\vk w}
\, d\vk{w},
\EQNY
Using  \kkk{change of variables} $\vk w =\squ \vk{b}_{j,l;u}-\vk{x}/ \squ$ we further obtain
\BQNY
p_{j,l;u}= u^{-1}\int_{\R^{2}}\phi_{ \SI_{ j,l;u} }( \squ \vk{b}_{j,l;u}-\vk{x}/ \squ ) P_{j,l;u}(\vk x) \, d\vk{x},
\EQNY
where
\BQNY
P_{j,l;u}(\vk x) :=\pk{ \underset{s\in[0,S]}{\exists_{t \in [0, T]}}
\begin{array}{ccc}
X_1(t_A+\frac{jT}{u}+\frac{t}{u})- X_1(t_A+\frac{jT}{u})>
 \frac{\mu_1}{\sqrt u} t+\frac{x_1}{\sqrt u}\\
X_2(s_A+\frac{lS}{u}+\frac{s}{u})- X_2(s_A+\frac{lS}{u})>  \frac{\mu_2}{\sqrt u} s +\frac{x_2}{\sqrt u}
\end{array}
\Bigg{|}  \vk Z_{j,l;u}= \squ \vk{b}_{j,l;u}-\frac{\vk{x}}{ \squ}}.
\EQNY
Now, we analyse $P_{j,l;u}(\vk x)$.
Due to the fact that $(t_A,s_A)\in A$, we have for all $t\in[0,T], s\in[0,S]$, and large enough $u$ 
\BQNY
t_A+\frac{jT}{u}+\frac{t}{u}\ge t_A+\frac{jT}{u}>s_A+\frac{lS}{u}+\frac{s}{u}\ge s_A+\frac{lS}{u}.
\EQNY
Thus, by the properties of Brownian motion  
\BQNY
P_{j,l;u}(\vk x) &=  &\pk{\exists_{t\in[0,T]}\  X_1(t)-\mu_1 t> x_1}\\
&&\times \pk{\exists_{s\in[0,S]} \ X_2(s_A+\frac{lS}{u}+\frac{s}{u})- X_2(s_A+\frac{lS}{u})> \frac{\mu_2}{\sqrt u} s + \frac{x_2}{\sqrt u} \Big{|} \vk Z_{j,l;u}= \squ \vk{b}_{j,l;u}-\frac{\vk{x}}{ \squ}},
\EQNY

Next we have
$$
\phi_{ \SI_{ j,l;u} }( \squ \vk{b}_{j,l;u}-\vk{x}/ \squ ) = \frac{1}{\sqrt{(2\pi)^{2}\abs{\SI_{j,l;u}}}}\exp\LT(-\frac{1}{2} ( \squ \vk{b}_{j,l;u}-\vk{x}/ \squ ) ^\top(\Sigma_{j,l;u})^{-1} ( \squ \vk{b}_{j,l;u}-\vk{x}/ \squ )   \RT),
 $$
where the exponent can be rewritten as
 \BQNY
 && ( \squ \vk{b}_{j,l;u}-\vk{x}/ \squ ) ^\top(\Sigma_{j,l;u})^{-1} ( \squ \vk{b}_{j,l;u}-\vk{x}/ \squ )  \\
 &&=u   (\vk{b}_{j,l;u})^\top\Sigma_{j,l;u}^{-1}  \vk{b}_{j,l;u}  - 2  \vk x^\top \Sigma_{j,l;u}^{-1}  \vk{b}_{j,l;u}   + \frac{1}{u } \vk{x}^\top \Sigma_{j,l;u}^{-1} \vk{x} \\
 &&=  u   g_A(t_A+\frac{jT}{u}, s_A+\frac{lS}{u})  - 2  \vk x^\top \Sigma_{j,l;u}^{-1}  \vk{b}_{j,l;u}   + \frac{1}{u } \vk{x}^\top \Sigma_{j,l;u}^{-1} \vk{x} .
 \EQNY
Define
\BQNY
f_{j,l;u}(\vk x) :=   \exp\LT( \vk x^\top\Sigma_{j,l;u}^{-1}  \vk{b}_{j,l;u}   - \frac{1}{2u } \vk{x}^\top\Sigma_{j,l;u}^{-1} \vk{x}    \RT),\ \ \ \ \vk x \in \R^2. 
\EQNY
Thus, it follows that
\BQNY 
p_1(u) =  \frac{u^{-1}}{ 2\pi}
\sum_{j=-N^{(1)}_u}^{N^{(1)}_u} \sum_{l=-N^{(2)}_u}^{N^{(2)}_u} \frac{1}{\sqrt{\abs{\SI_{j,l;u} }}} \exp\LT(-\frac{1}{2} u   g_A(t_A+\frac{jT}{u}, s_A+\frac{lS}{u})  \RT)
 \int_{\R^{2}}  f_{j,l;u}(\vk x)P_{j,l;u}(\vk x)\,d\vk{x}.\nonumber
\EQNY
Further, we obtain
from (i) of Lemma \ref{Lem:Taylor} that, for all large enough $u$,
\BQNY
g_A(t_A+\frac{jT}{u}, s_A+\frac{lS}{u}) \sim g_A(t_A,s_A)+\frac{1}{2}\LT(a_1 \LT(\frac{jT}{u}\RT)^2-2a_2  \LT(\frac{jT}{u}\RT)\LT(\frac{lS}{u}\RT) +a_3 \LT(\frac{lS}{u}\RT)^2\RT)
\EQNY
holds uniformly for   $-N_u^{(1)}\le j\le N_u^{(1)}, -N_u^{(2)}\le l\le N_u^{(2)}.$

 Consequently, by Lemma \ref{Lem:fPH}  below we  \kkk{obtain}
\BQNY 
\limit{u} \frac{p_1(u)}{ \exp\LT(-  g_A(t_A, s_A)u  /{2}    \RT) } =  \frac{1}{2\pi\sqrt{\abs{\SI_{0} }}} \frac{\H (\mu_1;T)\H(\mu_2-2\mu_1\rho;S)}{TS}
\int_{\R^2} e^{\frac{- (a_1x_1^2-2a_2 x_1x_2+a_3 x_2^2)}{4}} \,d\vk x,
\EQNY
 which gives the  result for $p_1(u)$. The claim for $p_2(u)$ follows with the same arguments.
\QED
\BEL\label{Lem:fPH}
\kkk{For any} $T,S>0$
\BQNY 
\lim_{u\to\IF}\int_{\R^{2}}  f_{j,l;u}(\vk x)P_{j,l;u}(\vk x)\,d\vk{x}=\H (\mu_1;T)\H(\mu_2-2\mu_1\rho;S)
\EQNY
holds uniformly for  $-N_u^{(1)}\le j\le N_u^{(1)}, -N_u^{(2)}\le l\le N_u^{(2)}.$
\EEL

\kkk{We omit the tedious proof of Lemma \ref{Lem:fPH} since its idea is standard,}
i.e., it is based on finding a uniform integrable bound for the integrand and then using the dominated convergence theorem.

\subsection{Proof of Lemma \ref{Lem:Pi12}}

Let us begin  with $\Pi_1(u)$.
 It follows that
\BQNY
\Pi_1(u)& = &\sum_{j=-N^{(1)}_u}^{N^{(1)}_u}\sum_{ -N^{(2)}_u\le l_1< l_2\le  N^{(2)}_u}p_{j,l_1,l_2;u}\\
&=&\sum_{j=-N^{(1)}_u}^{N^{(1)}_u} \sum_{l=-N^{(2)}_u}^{N^{(2)}_u}  p_{j,l,l+1;u}+ \sum_{j=-N^{(1)}_u}^{N^{(1)}_u} \sum_{l_1=-N^{(2)}_u}^{N^{(2)}_u} \sum_{  l_2=l_1+2} ^{N^{(2)}_u} p_{j,l_1,l_2;u} =: \Pi_{11}(u)+\Pi_{12}(u).
\EQNY
In order to deal with $\Pi_{11}(u)$ we note that
\BQNY
p_{j,l,l+1;u}=p_{j,l;u}+p_{j,l+1;u}-\widetilde p_{j,l;u},
\EQNY
where
\begin{equation*} %
\widetilde p_{j,l;u}=\pk{\exists_{ (t,s)\in\Del^{(1)}_{j;u}\times (\Del_{l;u}^{(2)}\cup \Del_{l+1;u}^{(2)})} \ X_1(t)>
\sqrt{u}(1+ \mu_1 t ), X_2(s)>
\sqrt{u}(1+ \mu_2 s )}.
\end{equation*}
Then we have
\BQNY
\Pi_{11}(u)=\sum_{j=-N^{(1)}_u}^{N^{(1)}_u} \sum_{l=-N^{(2)}_u}^{N^{(2)}_u} (p_{j,l;u}+p_{j,l+1;u}-\widetilde p_{j,l;u}).
\EQNY
\jj{Using the same arguments as in the proof of Lemma 4.2 we obtain}
\BQNY
\lim_{u\to\IF}\frac{\Pi_{11}(u)}{  e^{ -  g_A(t_A,s_A) u /2    } }=  \frac{1}{\mu_1 (\mu_2-2\mu_1\rho)}  \LT(\frac{2\H(\mu_1;T)\H(\mu_2-2\mu_1\rho;S)}{T S}  -\frac{\H(\mu_1;T)\H(\mu_2-2\mu_1\rho;2S)}{T S} \RT),
\EQNY
which gives that
\BQNY
\limsup_{S\to\IF}\limsup_{T\to\IF}\lim_{u\to\IF}\frac{\Pi_{11}(u)}{  e^{ - g_A(t_A,s_A) u /2  } }= 0.
\EQNY

Next we consider  $\Pi_{12}(u)$ which is  more involved. We have (recall \eqref{eq:bba} for $a_{j;u},b_{l;u}$)
\begin{eqnarray} 
p_{j,l_1,l_2;u}&=&\pk{ \underset{ s_2\in[0,S]} {\underset{ s_1  \in [0, S]} {\underset{ t  \in [0, T]} \exists}}
\begin{array}{ccc}
 X_1(t_A+\frac{jT}{u}+\frac{t}{u})> a_{j;u}
\sqrt{u} +\frac{\mu_1}{\sqrt u} t \\
 X_2(s_A+\frac{l_1S}{u}+\frac{s_1}{u})> \LLc{b_{l_1;u}}
\sqrt{u} +\frac{\mu_2}{\sqrt u} s_1\\
 X_2(s_A+\frac{l_2S}{u}+\frac{s_2}{u})> b_{l_2;u}
\sqrt{u} +\frac{\mu_2}{\sqrt u} s_2
\end{array}
} \nonumber\\
&\le &
\pk{ \underset{ s_2\in[0,S]} {\underset{ s_1  \in [0, S]} {\underset{ t  \in [0, T]} \exists}}
\begin{array}{cc}
X_1(t_A+\frac{jT}{u}+\frac{t}{u})> a_{j;u}
\sqrt{u} +\frac{\mu_1}{\sqrt u} t \\
\frac{1}{2}\LT(X_2(s_A+\frac{l_1S}{u}+\frac{s_1}{u})+ X_2(s_A+\frac{l_2S}{u}+\frac{s_2}{u})\RT)> b_{l_1,l_2;u}
\sqrt{u} +\frac{\mu_2}{2\sqrt u} (s_1+s_2)
\end{array}
}=: P_{j,l_1,l_2;u},
\end{eqnarray}
with
\BQNY
 b_{l_1,l_2;u} 
 =1+\mu_2\LT(s_A+\frac{l_1S}{u}+\frac{(l_2-l_1)S}{2u}\RT).
\EQNY
For notational simplicity, we shall denote
\BQNY
\widetilde{t_A}=t_A+\frac{jT}{u}, \ \   \widetilde{s_A}=s_A+\frac{l_1S}{u},\ \ \overline{s_A}=\widetilde{s_A}+\frac{(l_2-l_1)S}{2 u}, \ \  \widehat{s_A}=\widetilde{s_A}+\frac{(l_2-l_1)S}{4 u}.
\EQNY

Again by conditioning on the event 
\BQNY
E_{j,l_1,l_2;u}(x_1,x_2):=\LT\{X_1(\widetilde{t_A})=a_{j;u} \sqrt{u}-\frac{x_1}{\sqrt u}, \ \ \ \frac{1}{2}\LT(X_2(\widetilde{s_A})+ X_2(s_A+\frac{l_2S}{u} )\RT)=b_{l_1,l_2;u}
\sqrt{u} -\frac{x_2}{ \sqrt u}
\RT\},
\EQNY
we have 
\BQNY
P_{j,l_1,l_2;u}= u^{-1}\int_{\R^{2}}\phi_{ \SI_{ j,l_1,l_2;u} }( \squ \vk{b}_{j,l_1,l_2;u}-\vk{x}/ \squ )
F(j,l_1,l_2;u,\vk x)
\, d\vk{x},
\EQNY
where
\BQNY
\SI_{ j,l_1,l_2;u} =
\left(
   \begin{array}{cc}
   \widetilde{t_A}  & \rho\  \overline{s_A}\\
     \rho\ \overline{s_A} &  \widehat{s_A}\\
   \end{array}
 \right), \ \ \ \vk{b}_{j,l_1,l_2;u}=(a_{j;u}, b_{l_1,l_2;u})^\top
\EQNY
and
\BQNY
F(j,l_1,l_2;u,\vk x) := \pk{\exists_{ t  \in [0, T]}\  X_1(t)-\mu_1 t>
 x_1  } \  \pk{ \underset{ s_2\in[0,S]} {\underset{ s_1  \in [0, S]} { \exists}}
Y_{j,l_1,l_2;u}(s_1,s_2)  >x_2\
\Bigg{|}  E_{j,l_1,l_2;u}(x_1,x_2) },
\EQNY
where
\BQNY
Y_{j,l_1,l_2;u}(s_1,s_2)=\frac{\sqrt u}{2}
 \LT(
X_2(\widetilde{s_A}+\frac{s_1}{u})- X_2(\widetilde{s_A})+
X_2(s_A+\frac{l_2S}{u}+\frac{s_2}{u})- X_2(s_A+\frac{l_2S}{u}) \RT)-\frac{\mu_2}{2}(s_1+s_2).
\EQNY
Similarly as in the proof of Lemma \ref{Lem:PL}, we obtain
\BQNY
\phi_{ \SI_{ j,l_1,l_2;u} }( \squ \vk{b}_{j,l_1,l_2;u}-\vk{x}/ \squ )  =  \frac{1}{\sqrt{(2\pi)^{2}\abs{ \SI_{j,l_1,l_2;u} }}}\exp\LT(-\frac{1}{2} u \  (\vk{b}_{j,l_1,l_2;u})^\top\Sigma_{j,l_1,l_2;u}^{-1}  \vk{b}_{j,l_1,l_2;u}   \RT)\  f_{j,l_1,l_2;u}(\vk x),
\EQNY where
\BQNY
f_{j,l_1,l_2;u}(\vk x):= \exp\LT(   \vk x^\top \Sigma_{j,l_1,l_2;u}^{-1}  \vk{b}_{j,l_1,l_2;u}  - \frac{1}{2u } \vk{x}^\top \Sigma_{j,l_1,l_2;u}^{-1} \vk{x} \RT).
\EQNY
Next, some elementary calculations give that
 \BQNY
 (\vk{b}_{j,l_1,l_2;u})^\top\Sigma_{j,l_1,l_2;u}^{-1}  \vk{b}_{j,l_1,l_2;u} =    g_A(t_A+\frac{jT}{u}, s_A+\frac{l_1S}{u}+\frac{(l_2-l_1)S}{2 u}) + \frac{ \widetilde{t_A} g_A( \widetilde{t_A}, \overline{s_A} ) -a_{j;u} ^2  }{4  (\widetilde{t_A} \widehat{s_A}-\rho^2{\overline{s_A}}^2)}\ \frac{(l_2-l_1)S}{u}.
 \EQNY
Further, note that
\BQNY
g_A(t_A+\frac{jT}{u}, s_A+\frac{l_1S}{u}+\frac{(l_2-l_1)S}{2 u})= g_A(t_A+\frac{jT}{u}, s_A+\frac{l_1S}{u})+\frac{\partial g_A(t,s)}{\partial s}  \mid _{(\widetilde{t_A}, \widetilde{s_A}+\theta_{l_1,l_2;u}\frac{(l_2-l_1)S}{2 u})}  \frac{(l_2-l_1)S}{2 u}
\EQNY
holds for some $\theta_{l_1,l_2;u}\in(0,1)$ and
\BQNY
\frac{\partial g_A(t,s)}{\partial t}  \mid _{(\widetilde{t_A}, \widetilde{s_A}+\theta_{l_1,l_2;u}\frac{(l_2-l_1)S}{2 u})} \ \ \to\ 0, \ \ \ \ u\to\IF
\EQNY
holds uniformly for $ j, l_1,l_2 $ (hereafter when we write $ j, l_1,l_2 $ we mean $-N^{(1)}_u\le j \le N^{(1)}_u, -N^{(2)}_u\le   l_1,l_2\le N^{(2)}_u$)).

Consequently
\BQN\label{eq:expQ}
\exp\LT(-\frac{1}{2} u \  (\vk{b}_{j,l_1,l_2;u})^\top\Sigma_{j,l_1,l_2;u}^{-1}  \vk{b}_{j,l_1,l_2;u}   \RT) \sim
\exp\LT(-\frac{1}{2} u \ g_A(t_A+\frac{jT}{u}, s_A+\frac{l_1S}{u}) \RT) e^{- Q_0 (l_2-l_1)S}
\EQN
holds uniformly for $j, l_1,l_2$ \kkk{as $u\to\IF$},
where (by (b).(i) of Lemma 9 of \cite{DJT19a} \LLc{or  Lemma \ref{Lem:Taylor}.(i) \kkk{with} $a_1>0$}) 
\BQNY
Q_0=\frac{t_A g_A(t_A,s_A) -(1+\mu_1 t_A)^2}{8(t_A s_A-\rho^2s_A^2)} >0.
\EQNY
Next, we consider the uniform, in $ j,   l_1,l_2$, limit of the following: 
\BQNY
\pk{ \underset{ s_2\in[0,S]} {\underset{ s_1  \in [0, S]} { \exists}}
Y_{j,l_1,l_2;u}(s_1,s_2)  >x_2\
\Bigg{|}  E_{j,l_1,l_2;u}(x_1,x_2) }, \ \ \ u\to\IF
\EQNY
For the conditional mean we can derive that
\BQNY
\E{  Y_{j,l_1,l_2;u}(s_1,s_2)
{ \mid}  E_{j,l_1,l_2;u}(x_1,x_2) }&=& -\frac{\mu_2}{2}(s_1+s_2) \\
&&+ \LT(\frac{\rho (s_1+s_2)}{2\sqrt u}, \frac{s_1}{4\sqrt u}\RT)\ \Sigma_{j,l_1,l_2;u}^{-1} ( \vk{b}_{j,l_1,l_2;u}-\vk{x}/ \squ ),
\EQNY
which further gives \kkk{that}
\BQNY
&&\E{  Y_{j,l_1,l_2;u}(s_1,s_2)
{ \mid}  E_{j,l_1,l_2;u}(x_1,x_2) }= -\frac{\mu_2}{2}(s_1+s_2) + \frac{2\rho a_{j;u} \widehat{s_A}-\rho  a_{j;u} \overline{s_A} -2\rho^2  b_{j,l_1,l_2;u} \overline{s_A} + b_{j,l_1,l_2;u} \widetilde{t_A}}{4(\widetilde{t_A} \widehat{s_A}-\rho^2 {\overline{s_A}}^2)} s_1\\
&& \ \ \ \  \ \ \ \ \ \ \ +\frac{ \rho a_{j;u} \widehat{s_A}-\rho^2  b_{j,l_1,l_2;u} \overline{s_A} }{2(\widetilde{t_A} \widehat{s_A}-\rho^2 {\overline{s_A}}^2)} s_2
+ \frac{ \rho  \overline{s_A} s_1 -2\rho   \widehat{s_A} (s_1+s_2)  }{4(\widetilde{t_A} \widehat{s_A}-\rho^2 {\overline{s_A}}^2)}\  \frac{x_1}{u}
+ \frac{ 2\rho ^2  \overline{s_A} (s_1+s_2) -  \widetilde{t_A}s_1 }{4(\widetilde{t_A} \widehat{s_A}-\rho^2 {\overline{s_A}}^2)}\  \frac{x_2}{u}\\
&& \ \ \ \  \ \ \ \ \ \ \ \to -\frac{1}{2}\LT(\mu_2-2\mu_1\rho\RT) s_2, \ \ \ \ u\to\IF.
\EQNY
For the conditional variance of the increments we have
\BQNY
&&\mathrm{Var}\LT\{  Y_{j,l_1,l_2;u}(s_1,s_2)  -  Y_{j,l_1,l_2;u}(s'_1,s'_2)
{ \mid}  E_{j,l_1,l_2;u}(x_1,x_2) \RT\} = \frac{\abs{s_1-s'_1}+\abs{s_2-s'_2}}{4}\\
&&\ \ \ \ \ \ \  \ \ \ \ \ \ \ \ \ \ \ \ \ +\LT(\frac{\rho(s_1-s'_1 + s_2-s'_2)}{2\sqrt u}, \frac{s_1-s'_1}{4 \sqrt u}\RT) \Sigma_{j,l_1,l_2;u}^{-1} \LT(\frac{\rho(s_1-s'_1 + s_2-s'_2)}{2\sqrt u}, \frac{s_1-s'_1}{4 \sqrt u}\RT)^{\top}\\
&&\ \ \ \ \ \ \  \ \ \ \ \ \ \ \ \ \ \ \ \ \to  \frac{\abs{s_1-s'_1}+\abs{s_2-s'_2}}{4}, \ \ \ u\to\IF.
\EQNY
Therefore, \LLc{similarly as in Lemma \ref{Lem:fPH} we can show that as $u\to\IF$}
\BQNY
\pk{ \underset{ s_2\in[0,S]} {\underset{ s_1  \in [0, S]} { \exists}}
Y_{j,l_1,l_2;u}(s_1,s_2)  >x_2\
\Bigg{|}  E_{j,l_1,l_2;u}(x_1,x_2) } \  \to\ \pk{\underset{ s_2\in[0,S]} {\underset{ s_1  \in [0, S]} { \exists}} \frac{1}{2}(B_1(s_1)+B_2(s_2)) -\frac{1}{2}\LT( \mu_2-2\mu_1\rho\RT) s_2 >x_2}.
\EQNY

\jj{Consequently, the dominated convergence theorem gives}
\BQN\label{eq:fF1}
 &&\int_{\R^{2}}  f_{j,l_1,l_2;u}(\vk x) F(j,l_1,l_2;u,\vk x) \,d\vk{x} \nonumber\\
 && \to
 \int_{\R}e^{2\mu_1 x_1} \pk{\exists_{t\in[0,T]}\ X_1(t)-\mu_1 t> x_1}dx_1 \nonumber\\
&&\ \ \ \ \ \times \int_\R e^{2(\mu_2-2\mu_1\rho) x_2} \pk{\underset{ s_2\in[0,S]} {\underset{ s_1  \in [0, S]} { \exists}} \frac{1}{2}(B_1(s_1)+B_2(s_2)) -\frac{1}{2}\LT(\mu_2-2\mu_1\rho\RT) s_2 >x_2}dx_2\\
&&=: \H(\mu_1; T)\ \H(\mu_1,\mu_2;S) \nonumber
\EQN
holds uniformly for $j,l_1,l_2$, as $u\to\IF$.

 \kkk{
Next we derive a useful upper bound for $\H(\mu_1,\mu_2;S)$, $S>0$:
\BQN
\H(\mu_1,\mu_2;S)   \le  (\lfloor S\rfloor)^2 e^{ Q_0 S} \H(\mu_1,\mu_2;1)<\IF.
\label{Lem:Hmumu}
\EQN
In order to  prove (\ref{Lem:Hmumu}), by taking
$j=l_1=0, l_2=1$ 
we arrive at
\BQN\label{eq:PU}
P_{0,0,1;u}&=&\pk{ \underset{ s_2\in[0,S]} {\underset{ s_1  \in [0, S]} {\underset{ t  \in [0, T]} \exists}}
\begin{array}{cc}
X_1(t_A+\frac{t}{u})> a_{0;u}
\sqrt{u} +\frac{\mu_1}{\sqrt u} t \\
\frac{1}{2}(X_2(s_A+\frac{s_1}{u})+ X_2(s_A+\frac{S}{u}+\frac{s_2}{u}))> b_{0,1;u}
\sqrt{u} +\frac{\mu_2}{2\sqrt u} (s_1+s_2)
\end{array}
}\nonumber \\
&\sim & \frac{u^{-1}}{\sqrt{(2\pi)^{2}\abs{ \SI_{0,0,0;u} }}} \exp\LT(-\frac{1}{2} u \ g_A(t_A, s_A) \RT) e^{- Q_0 S} \H(\mu_1; T)\H(\mu_1,\mu_2;S).
\EQN
Define, for any integers $0\le m,n\le \lfloor S\rfloor$,
\BQNY
q_{m,n;u} := \pk{ \underset{ s_2\in[0,1]} {\underset{ s_1  \in [0, 1]} {\underset{ t  \in [0, T]} \exists}}
\begin{array}{cc}
X_1(t_A+\frac{t}{u})> a_{0;u}
\sqrt{u} +\frac{\mu_1}{\sqrt u} t \\
\frac{1}{2}(X_2(s_A+\frac{m}{u}+\frac{s_1}{u})+ X_2(s_A+\frac{S+n}{u}+\frac{s_2}{u}))> \widetilde b_{m,n;u}
\sqrt{u} +\frac{\mu_2}{2\sqrt u} (s_1+s_2)
\end{array}
}
\EQNY
with
\BQNY
\widetilde b_{m,n;u}=1+\mu_2\LT(s_0+\frac{m}{u}+\frac{ S+n-m}{2u}\RT).
\EQNY
Using the same arguments as in the derivation of \eqref{eq:PU} one can show that
\BQN\label{eq:qq}
q_{m,n;u} \sim \frac{u^{-1}}{\sqrt{(2\pi)^{2}\abs{ \SI_{0,0,0;u} }}} \exp\LT(-\frac{1}{2} u \ g_A(t_A, s_A) \RT) e^{- Q_0 (S+n-m)} \H(\mu_1; T)\H(\mu_1,\mu_2;1).
\EQN
Comparing \eqref{eq:PU} and \eqref{eq:qq} we derive
\BQNY 
\H(\mu_1,\mu_2;S)  &\le& \sum_{m=0}^{\lfloor S\rfloor-1} \sum_{n=0}^{\lfloor S\rfloor-1} e^{- Q_0 (n-m)} \H(\mu_1,\mu_2;1)\nonumber\\
&\le & (\lfloor S\rfloor)^2 e^{ Q_0 S} \H(\mu_1,\mu_2;1).
\EQNY
The finiteness of $\H(\mu_1,\mu_2;1)$ can be proved by  using the Borell-TIS inequality.
This justifies
bound (\ref{Lem:Hmumu}).
}

Now, we are ready to analyse the triple sum $\Pi_{12}(u)$.
We have
\BQNY
\Pi_{12}(u)& = &\sum_{j=-N^{(1)}_u}^{N^{(1)}_u} \sum_{l_1=-N^{(2)}_u}^{N^{(2)}_u} \sum_{  l_2=l_1+2} ^{N^{(2)}_u}\frac{u^{-1}}{\sqrt{(2\pi)^{2}\abs{\SI_{j,l_1,l_2;u} }}} \\
&& \ \ \ \ \ \ \ \times\exp\LT(-\frac{1}{2} u \  (\vk{b}_{j,l_1,l_2;u})^\top\Sigma_{j,l_1,l_2;u}^{-1}  \vk{b}_{j,l_1,l_2;u}   \RT)
 \int_{\R^{2}}  f_{j,l_1;u}(\vk x) F(j,l_1,l_2;u,\vk x) \,d\vk{x}. 
\EQNY
Therefore, we can derive from   \eqref{eq:expQ}-\eqref{eq:fF1}
\kkk{and (\ref{Lem:Hmumu})}
that
\BQNY
\lim_{u\to\IF}\frac{ \Pi_{12}(u)}{\exp(-ug_A(t_A,s_A)/2)} \le \mathrm{Const}  \  \sum_{ k=1} ^{\IF}  e^{- kQ_0 S} \frac{{\H(\mu_1;T) \H(\mu_1,\mu_2;1)(\lfloor S\rfloor)^2}}{TS} .
\EQNY
Consequently, the above implies that
\BQNY
\limsup_{S\to\IF} \limsup_{T\to\IF}\lim_{u\to\IF}\frac{\Pi_{12}(u)}{\exp(-ug_A(t_A,s_A)/2)}=0.
\EQNY
Thus, the claim for $\Pi_{1}(u)$ is established. \LLc{Using similar arguments, one can further show that the claim for $\Pi_2(u)$ holds. }\QED


\subsection{Proof of Lemma \ref{Lem:r123}}

The claim for $r_1(u)$ follows from the same arguments as that for $r_0(u)$ of \nelem{Lem:r01}.
Next, as in the proof of \nelem{Lem:r01}, using the Piterbarg's inequality we can show that
\BQNY 
r_2(u)\le    C_2 u^{3/2} e^{-\frac{  u }{2} \widetilde g_u},
\EQNY
where $C_2>0$ is some constant which does not depend on $u,$ and \LLc{thus the claim for $r_2(u)$ follows since}
\BQNY
\widetilde g_u&=& \inf_{(t,s)\in D_1\cup D_2 } g(t,s)=\inf_{s \in[s_0-\theta_0, s_0-\ln(u)/\sqrt u]\cup[ s_0+\ln(u)/\sqrt u, s_0+\theta_0]} g_L(s)\\
&\ge&
  g_L(s^*) +\frac{b_0}{2}(1-\vn) \frac{(\ln(u))^2}{u},
\EQNY
where the last inequality follows by (ii.3) of Lemma \ref{Lem:Taylor}.
\kkk{Finally}, the claim for  $r_3(u)$ can be proved similarly, by using Piterbarg's inequality and (ii.1)-(ii.2) of Lemma \ref{Lem:Taylor}. \QED

\subsection{Proof of Lemma \ref{Lem:PL2}}

We first analyse the summand $p_{j;u}$. \kkk{Let} 
\BQNY
\vk b_{j;u}=(a_{j;u}, b_{j,u})^\top, \ \ \   a_{j;u}=1+\mu_1 (t^*+\frac{jT}{u}), \ \ \
b_{j;u}=1+\mu_2 (s^*+\frac{jT}{u}).
\EQNY
\kkk{Then} 
\begin{eqnarray*} 
p_{j;u}
 = \pk{  \underset{ (t,s)  \in \Del_{T,S} }\exists
\begin{array}{ccc}
 X_1(t^*+\frac{jT}{u}+\frac{t}{u})> a_{j;u}
\sqrt{u} +\frac{\mu_1}{\sqrt u} t \\
 X_2(s^*+\frac{jT}{u}+\frac{s}{u})> b_{j;u}
\sqrt{u} +\frac{\mu_2}{\sqrt u} s
\end{array}
}. 
\end{eqnarray*}
 Define $\vk Z_{j;u}:=(X_1(t^*+\frac{jT}{u}),  X_2(s^*+\frac{jT}{u}))^\top$, whose  density function   is given by
\BQNY
\phi_{\SI_{j;u}}( \vk w)=\frac{1}{\sqrt{(2\pi)^{2}\abs{ \SI_{j;u} }}}\exp\LT(-\frac{1}{2}\vk w^\top(\Sigma_{j;u})^{-1}  \vk w  \RT),\ \ \ \vk w=(w_1,w_2)^\top,
\EQNY
with
the covariance matrix  given by
$$
\  \Sigma_{j;u}=
\left(
   \begin{array}{cc}
     t^*+\frac{jT}{u}  & \rho\  (t^*+\frac{jT}{u})\\
     \rho\  (t^*+\frac{jT}{u}) & s^*+\frac{jT}{u} \\
   \end{array}
 \right).
$$

By conditioning on the value of $\vk Z_{j;u}$ and using  \kkk{change of variables} $\vk w =\squ \vk{b}_{j;u}-\vk{x}/ \squ$, we further obtain
\BQNY
p_{j;u}= u^{-1}\int_{\R^{2}}\phi_{ \SI_{ j;u} }( \squ \vk{b}_{j;u}-\vk{x}/ \squ )
\pk{   \underset{ (t,s)  \in \Del_{T,S} }\exists
\begin{array}{ccc}
X_1(t)- \mu_1 t>x_1\\
X_2(s)- \mu_2 s>x_2
\end{array}
}
\, d\vk{x}.
\EQNY
Consequently, similar arguments as in the proof of Lemma \ref{Lem:PL} yield
\BQNY 
p_1(u) &\sim& p_2(u)\ \sim\  \sum_{j=-N_u^{(1)}}^{N_u^{(1)}} p_{j;u}
 \sim   \frac{\H(T,S)u^{-1}}{\sqrt{(2\pi t^*)^2(1-\rho^2)}}  \sum_{j=-N_u^{(1)}}^{N_u^{(1)}} e^{-\frac{u}{2} g_L(t^*+\frac{jT}{u})}\nonumber\\
 &\sim&     \frac{\H(T,S)u^{-1/2}}{T\sqrt{(2\pi t^*)^2(1-\rho^2)}}   e^{-\frac{u}{2} g_L(t^*)}\int_{\R}e^{-\frac{b_0}{4}x^2} dx.
\EQNY
This completes the proof. \QED

\subsection{Proof of Lemma \ref{Lem:subadd}}
First note that
\BQNY
&&\pk{\exists_{ ( t,s) \in (t^*,s^*)+u^{-1}\Del_{T_1+T_2,S}} \ X_1(t)>
\sqrt{u}(1+ \mu_1 t ), X_2(s)>
\sqrt{u}(1+ \mu_2 s )}\nonumber\\
&&\ \ \le \pk{\exists_{ ( t,s) \in (t^*,s^*)+u^{-1}\Del_{T_1,S}} \ X_1(t)>
\sqrt{u}(1+ \mu_1 t ), X_2(s)>
\sqrt{u}(1+ \mu_2 s )}\nonumber\\
&&\ \ \ \ + \pk{\exists_{ ( t,s) \in (t^*+\frac{T_1}{u},s^*+\frac{T_1}{u})+u^{-1}\Del_{T_2,S}} \ X_1(t)>
\sqrt{u}(1+ \mu_1 t ), X_2(s)>
\sqrt{u}(1+ \mu_2 s )}.
\EQNY
Using the same arguments as the proof of Lemma \ref{Lem:PL2}, we conclude the sub-additivity of $\H(T,S), T>0.$
\kkk{Thus}
\BQNY
\lim_{T\to\IF}\frac{1}{T}\H(T,S) = \inf_{T> 0}\frac{1}{T}\H(T,S) <\IF
\EQNY
  follows directly from Fekete's lemma. Moreover, since by definition
\BQNY
\H(T,S)\ge \int_{\R^2}e^{\vk x^\top \Sigma_{*}^{-1}  \vk{b}_{*} }\pk{  {\underset{  t \in [0,T]}\exists }
\begin{array}{ccc}
X_1(t)- \mu_1  t >  x_1 \\
X_2(t)-  \mu_2  t  >x_2
\end{array}
 } dx_1 dx_2,
\EQNY
  the positive lower bound follows from Lemma 4.7 in \cite{DHJT18}. This completes the proof. \QED

\subsection{Proof of Lemma \ref{Lem:Pi112}}
We \kkk{begin} with the analysis of $\Pi_1(u)$.
We first look at $p_{j,l;u}$. Denote
\BQNY
&&\vk b_{u}=\vk b_{j,l,m,n;u}:=(a_{j,m;u}, b_{j,l,m,n;u})^\top, \\
&& a_{j,m;u}=1+\mu_1 (t^*+\frac{jT+m}{u}), \ \ b_{j,l,m,n;u}=1+\mu_2 (t^*+\frac{jT+m}{u}+\frac{lS+n}{u}).
\EQNY
It is derived that
\BQNY
p_{j,l;u} 
&\le &\sum_{m=0}^{\lfloor T \rfloor-1}\sum_{n=0}^{\lfloor S \rfloor-1}
\pk{  \underset{s-t\in[\frac{lS}{u}+\frac{n}{u}, \frac{lS}{u}+\frac{n+1}{u} ]}{\underset{ t  \in [t^*+\frac{jT}{u}+\frac{m}{u}, t^*+\frac{jT}{u}+\frac{m+1}{u}] }\exists }
\ X_1(t)>
\sqrt{u}(1+ \mu_1 t ), X_2(s)>
\sqrt{u}(1+ \mu_2 s )}\\
&=& \sum_{m=0}^{\lfloor T \rfloor-1}\sum_{n=0}^{\lfloor S \rfloor-1}
\pk{  \underset{s\in[0,1]}{\underset{ t  \in [0, 1] }\exists }
\begin{array}{ccc}
X_1(t^*+\frac{jT+m}{u}+\frac{t}{u}) > a_{j,m;u}
\sqrt{u} +\frac{\mu_1}{\sqrt u} t \\
X_2(t^*+\frac{jT+m}{u}+\frac{lS+n}{u}+\frac{t+s}{u}) > b_{j,l,m,n;u}
\sqrt{u} +\frac{\mu_2}{\sqrt u} (t+s)
\end{array}
}\\
&=:& \sum_{m=0}^{\lfloor T \rfloor-1}\sum_{n=0}^{\lfloor S \rfloor-1} p_{j,l,m,n;u} .
\EQNY

Next, we look at $p_{j,l,m,n;u}$.
We define
\BQNY
\vk Z_{u}
&:=& \LT(X_1(t^*+\frac{jT+m}{u}),  \ X_2(s^*+\frac{jT+m}{u}+\frac{lS+n}{u})\RT)^\top,\\
Y_{1;u}(t)
&:= &\LT(X_1(t^*+\frac{jT+m}{u}+\frac{t}{u})- X_1(t^*+\frac{jT+m}{u})\RT) {\sqrt u} -\mu_1 t, \\
Y_{2;u}(t,s)
&:=&\LT(X_2(t^*+\frac{jT+m}{u}+\frac{lS+n}{u}+\frac{t+s}{u})- X_2(t^*+\frac{jT+m}{u}+\frac{lS+n}{u})\RT) {\sqrt u} -\mu_2 (t+s).
\EQNY
Consider the conditional process
\BQNY
\vk W_{u}(t,s):=(Y_{1;u}(t), Y_{2;u}(t,s))^\top\mid \vk Z_{u}= \squ \vk{b}_{u}-\frac{\vk{x}}{ \squ}.
\EQNY
We have that
$ (Y_{1;u}(t), Y_{2;u}(t,s), \vk Z_{u}^\top)$ is a normally distributed random vector,
with mean
\BQNY
\vk {\widehat\mu}(t,s) := (-\mu_1 t, -\mu_2(t+s), 0, 0)^\top
\EQNY
and  covariance matrix  given by (suppose $S>1$)
\BQNY
\widehat\Sigma_{u}(t,s):=
\left(
   \begin{array}{cccr}
   t&0&  0 &  \rho\frac{t}{\sqrt u}\\
0 & t+s & 0 & 0\\
0  & 0&  t^*+\frac{jT+m}{u}  & \rho\  (t^*+\frac{jT+m}{u})\\
\rho\frac{t}{\sqrt u} & 0&    \rho\  (t^*+\frac{jT+m}{u}) & t^*+\frac{jT+m}{u} +\frac{lS+n}{u}\\
   \end{array}
 \right).
\EQNY
Thus, for the mean
\BQNY
\E{\vk W_{u}(t,s)}&=&  (-\mu_1 t, -\mu_2(t+s))+
\left(
   \begin{array}{cc}
   0  &  \rho\frac{t}{\sqrt u}\\
0&0\\
   \end{array}
 \right)
 \left(
   \begin{array}{cc}
     t^*+\frac{jT+m}{u}  & \rho\  (t^*+\frac{jT+m}{u})\\
     \rho\  (t^*+\frac{jT+m}{u}) & t^*+\frac{jT+m}{u}+\frac{lS+n}{u} \\
   \end{array}
 \right)^{-1} \LT( \squ \vk{b}_{u}-\frac{\vk{x}}{ \squ}\RT)\\
 &=& (-\mu_1 t, -\mu_2(t+s))+ \LT(\frac{\rho (b_{j,l,m,n;u} - \rho a_{j,m;u} ) t -  \rho   t \frac{x_2-\rho x_1}{u}}{(t^*+\frac{jT+m}{u}+\frac{lS+n}{u}) -\rho^2(t^*+\frac{jT+m}{u})} , 0\RT) \\ 
 &\to &  \LT(- \nu(\rho) \ t, -\mu_2   (t+s) \RT)^\top, \ \ \nu(\rho):=\frac{\rho^2-\rho+(\mu_1-\mu_2\rho)t^*}{t^*(1-\rho^2)},
\EQNY
as $u\to\IF$, where the convergence is uniform for $-N_u^{(1)}\le j\le N_u^{(1)}, 1\le l\le N_u^{(2)}$.
Similarly, we can derive that, for any $t_1,t_2\in[0,T], s_1,s_2\in[-S,S],$
\BQNY
\mathrm{Cov}(\vk W_{u}(t_1,s_1) -\vk W_{u}(t_2,s_2))
&\to &\mathrm{Cov}((B_1(t_1)-B_1(t_2), B_2(t_1+s_1)-B_2(t_2+s_2))^\top)
\EQNY
as $u\to\IF$, uniformly for $-N_u^{(1)}\le j\le N_u^{(1)}, 1\le l\le N_u^{(2)}$. Consequently,  we have, as $u\to\IF,$
\BQNY
&&\pk{  \underset{s\in[0,1]}{\underset{ t  \in [0, 1] }\exists }
\begin{array}{ccc}
Y_{1;u}(t)>x_1 \\
Y_{2;u}(t,s) >  x_2
\end{array}
 \Bigg{|} \vk Z_{u}= \squ \vk{b}_{u}-\frac{\vk{x}}{ \squ}}\\
&& \to \ \pk{  \underset{s\in[0,1]}{\underset{ t  \in [0, 1] }\exists }
\begin{array}{ccc}
B_1(t)- \nu(\rho)  t >  x_1 \\
B_2(t+s) - \mu_2 (t+s) >x_2
\end{array}
 }.
\EQNY
Similar arguments as those in the proof of Lemma \ref{Lem:PL} gives that
\BQNY
p_{j,l,m,n;u}  \sim   \frac{\widehat\H(1,1)u^{-1}}{\sqrt{(2\pi)^2\abs{\Sigma_*}}}   e^{-\frac{u}{2} g_B(t^*+\frac{jT+m}{u}, s^*+\frac{jT+m}{u}+\frac{lS+n}{u})},
\EQNY
where (recall notation in \eqref{eq:Sig-b})
\BQNY
\widehat\H(1,1)=\int_{\R^2} e^{\vk x^\top \Sigma_{*}^{-1}  \vk{b}_{*} } \pk{  \underset{s\in[0,1]}{\underset{ t  \in [0, 1] }\exists }
\begin{array}{ccc}
B_1(t)- \nu(\rho) t >  x_1 \\
B_2(t+s)  - \mu (t+s) >x_2
\end{array}
 } dx_1 dx_2\in(0,\IF).
\EQNY
It follows further from (ii.2) of Lemma \ref{Lem:Taylor} that  there exists some $\vn>0$ such that, for all $t<s$ small,
\BQNY
g_B(t^*+ t,t^*+ s)
 \ge  g_L(t^*)+b_2(1-\vn)( s-t )+\frac{c_2 (1-\vn) }{2} t^2,
\EQNY
thus, for $u$ sufficiently large
\BQNY
g_B(t^*+\frac{jT+m}{u}, t^*+\frac{jT+m}{u}+\frac{lS+n}{u}) \ge  g_L(t^*)+b_2(1-\vn)\frac{lS}{u}+\frac{c_2(1-\vn)}{2} \LT(\frac{\widehat j T}{u}\RT)^2
\EQNY
holds for all $-N_u^{(1)}\le j\le N_u^{(1)}, 1\le l\le N_u^{(2)}, 0\le m\le \lfloor T\rfloor-1, 0\le n\le \lfloor S\rfloor-1$, where $\widehat j=j I_{(j\ge 0)} +(j+1) I_{(j<0)}$. This implies that, for $u$ large
\BQNY
 e^{-\frac{u}{2} g_A(t^*+\frac{jT+m}{u}, s^*+\frac{jT+m}{u}+\frac{lS+n}{u})}\le e^{-\frac{u}{2} g_L(t^*)} \frac{\sqrt u}{T} \LT(   e^{- \frac{c_2(1-\vn)}{4} \LT(\frac{\widehat jT}{\sqrt u}\RT)^2} \frac{T}{\sqrt u}\RT)   e^{-\frac{b_2(1-\vn)}{2} lS}.
\EQNY
Based on the above discussions we obtain
\BQNY
\lim_{u\to\IF}\frac{\Pi_1(u)}{u^{-1/2} \exp(-  g_L(t^*){u}/{2})  } \le
    \frac{\widehat\H(1,1)}{\sqrt{(2\pi)^2\abs{\Sigma_*}}}    \ \frac{\lfloor T  \rfloor \lfloor S  \rfloor}{T}   \sum_{ l\ge 1 } e^{-\frac{b_2(1-\vn)}{2} lS} \int_{\R} e^{-\frac{c_2(1-\vn)}{4}   x^2}dx.
\EQNY
Similar bounds can be found for $\overline \Pi_1(u)$, and thus the first claim follows.

\COM{
\underline{Analysis of $\Pi_2(u) $}. Recall

\BQNY
q_{j_1,j_2;u}=\pk{
\begin{array}{cc}
\exists_{  ( t,s) \in (t_0+\frac{j_1 T}{u},t_0+\frac{j_1 T}{u})+u^{-1}\Del_{T,S}} \ X_1(t)>
\sqrt{u}(1+ \mu_1 t ), X_2(s)>
\sqrt{u}(1+ \mu_2 s ) \\
 \exists_{ ( t,s) \in (t_0+\frac{j_2 T}{u},t_0+\frac{j_2 T}{u})+u^{-1}\Del_{T,S}} \ X_1(t)>
\sqrt{u}(1+ \mu_1 t ), X_2(s)>
\sqrt{u}(1+ \mu_2 s )
 \end{array}
 }
\EQNY
and
\BQNY
\Pi_2(u) &=&  \sum_{ -N_u^{(1)}\le j_1< j_2\le  N_u^{(1)}}  q_{j_1,j_2;u} \\
&=& \sum^{ N_u^{(1)}} _{ j_1=-N_u^{(1)}} \sum_{j_2>j_1+1} q_{j_1,j_2;u} + \sum^{ N_u^{(1)}} _{ j_1=- N_u^{(1)}}   q_{j_1,j_1+1;u}=: \Pi_{21}(u)+\Pi_{22}(u)
\EQNY}
Next we consider $\Pi_{21}(u)$. For any $j_2>j_1+1$ we have
\BQNY
q_{j_1,j_2;u}&=&\pk{
\begin{array}{cc}
\exists_{  ( t,s) \in (t^*+\frac{j_1 T}{u},s^*+\frac{j_1 T}{u})+u^{-1}\Del_{T,S}} \ X_1(t)>
\sqrt{u}(1+ \mu_1 t ), X_2(s)>
\sqrt{u}(1+ \mu_2 s ) \\
 \exists_{ ( t,s) \in (t^*+\frac{j_2 T}{u},s^*+\frac{j_2 T}{u})+u^{-1}\Del_{T,S}} \ X_1(t)>
\sqrt{u}(1+ \mu_1 t ), X_2(s)>
\sqrt{u}(1+ \mu_2 s )
 \end{array}
 }\\
 &\le &
u^{-1}\int_{\R^{2}}\phi_{ \SI_{j_1,j_2;u} }( \squ \vk{b}_{j_1,j_2;u}-\vk{x}/ \squ )
\overline F(j_1,j_2;u,\vk x)
\, d\vk{x}=:Q_{j_1,j_2;u},
\EQNY
where, with $a_{j;u}=1+\mu_1(t^*+\frac{j_1T}{u}), b_{j;u}=1+\mu_2(s^*+\frac{j_1T}{u})$,
\BQNY
\SI_{j_1,j_2;u} =  \LT(t^*+\frac{j_1T}{u}+\frac{(j_2-j_1)S}{4u} \RT)
\left(
   \begin{array}{cc}
    1 & \rho \\
     \rho  & 1 \\
   \end{array}
 \right), \ \
\vk{b}_{j_1,j_2;u} =\LT(\frac{a_{j_1;u}+a_{j_2;u}}{2}, \frac{b_{j_1;u}+b_{j_2;u}}{2}\RT),
\EQNY
\BQNY
\overline F(j_1,j_2;u,\vk x):= \pk{  \underset{(t',s')\in\Del_{T,S}}{\underset{ (t,s)\in\Del_{T,S} }\exists }
\begin{array}{ccc}
Y_{1;u}(t,t') >  x_1 \\
Y_{2;u}(s,s') >x_2
\end{array}
\Bigg|
\begin{array}{ccc}
Y_{3;u} = \frac{a_{j_1;u}+a_{j_2;u}}{2}\sqrt u -\frac{ x_1}{\sqrt u} \\
Y_{4;u} =\frac{a_{j_1;u}+a_{j_2;u}}{2}\sqrt u -\frac{ x_2}{\sqrt u}
\end{array}
 },
\EQNY
with
\BQNY
&&Y_{1;u}(t,t')= \frac{\sqrt u}{2} \LT(X_1(t^*+\frac{j_1T}{u}+\frac{t}{u}) -  X_1(t^*+\frac{j_1T}{u})  + X_1(t^*+\frac{j_2T}{u}+\frac{t'}{u})-  X_1(t^*+\frac{j_2T}{u})\RT) -\frac{\mu_1}{2}(t+t'), \\
&&Y_{2;u}(s,s')= \frac{\sqrt u}{2} \LT(X_2(s^*+\frac{j_1T}{u}+\frac{s}{u})-  X_2(s^*+\frac{j_1T}{u})  + X_2(s^*+\frac{j_2T}{u}+\frac{s'}{u})-  X_2(s^*+\frac{j_2T}{u})\RT) -\frac{\mu_2}{2}(s+s'), \\
&&Y_{3;u} =\frac{1}{2} \LT(X_1(t^*+\frac{j_1T}{u})+ X_1(t^*+\frac{j_2T}{u})\RT), \ \  \ \
Y_{4;u}  = \frac{1}{2} \LT(X_2(s^*+\frac{j_1T}{u})+X_2(s^*+\frac{j_2T}{u})\RT).
\EQNY
Next we have that $(Y_{1;u}(t,t'), Y_{2;u}(s,s'), Y_{3;u},Y_{4;u})$  is a normally distributed random vector,
with mean
\BQNY
\vk {\widehat\mu}(t,t',s,s') = \LT(-\frac{\mu_1}{2}(t+t'), -\frac{\mu_2}{2}(s+s'), 0, 0\RT)^\top
\EQNY
and  covariance matrix  given by (suppose $ T>S$)
\BQNY
\widehat\Sigma_{u}(t,s)=
\left(
   \begin{array}{cccr}
  \frac{t+t'}{4}&\frac{\rho(t\wedge s+t' \wedge s')}{4}&   \frac{t}{4 \sqrt u}  &  \frac{\rho t}{4 \sqrt u}\\
\frac{\rho(t\wedge s+t' \wedge s')}{4} &  \frac{s+s'}{4}& \frac{\rho s}{4 \sqrt u} &  \frac{s}{4 \sqrt u}\\
\frac{t}{4 \sqrt u}  &  \frac{\rho s}{4 \sqrt u}&  t^*+\frac{j_1T}{u}+\frac{(j_2-j_1)T}{4u} & \rho\  \LT(t^*+\frac{j_1T}{u}+\frac{(j_2-j_1)T}{4u}\RT)\\
\frac{\rho t}{4 \sqrt u}  &  \frac{s}{4 \sqrt u}&    \rho\  \LT(t^*+\frac{j_1T}{u}+\frac{(j_2-j_1)T}{4u}\RT)   & s^*+\frac{j_1T}{u}+\frac{(j_2-j_1)T}{4u} \\
   \end{array}
 \right).
\EQNY
\COM{
Thus, for the mean function
\BQNY
\E{\vk W_{u}(t,t',s,s')}&=&  (-\mu_1 t, -\mu(t+s))+
\left(
   \begin{array}{cc}
   0  &  \rho\frac{t}{\sqrt u}\\
0&0\\
   \end{array}
 \right)
 \left(
   \begin{array}{cc}
     t_0+\frac{jT+m}{u}  & \rho\  (t_0+\frac{jT+m}{u})\\
     \rho\  (t_0+\frac{jT+m}{u}) & t_0+\frac{jT+m}{u}+\frac{lS+n}{u} \\
   \end{array}
 \right)^{-1} \LT( \squ \vk{b}_{u}-\frac{\vk{x}}{ \squ}\RT)\\
 &=&  \\ 
 &\to &  \LT(-\frac{1-\rho}{1+\rho}\mu_1 t, -\mu   (t+s) \RT)^\top,
\EQNY
where the convergence is uniform for $-N_u^{(1)}\le j\le N_u^{(1)}, ...$.
Similarly, we can derive that, for any $t_1,t_2\in[0,T], s_1,s_2\in[-S,S],$
\BQNY
\mathrm{Cov}(\vk W_{u}(t_1,s_1) -\vk W_{u}(t_2,s_2))
&\to &\mathrm{Cov}((X_1(t_1)-X_1(t_2), X_2(t_1+s_1)-X_2(t_2+s_2))^\top)
\EQNY
uniformly for $-N_u^{(1)}\le j\le N_u^{(1)}, ...$. Consequently,  we have 
\BQNY
&&\pk{  \underset{s\in[0,1]}{\underset{ t  \in [0, 1] }\exists }
\begin{array}{ccc}
Y_{1;u}(t)>x_1 \\
Y_{2;u}(t,s) >  x_2
\end{array}
 \mid  \vk Z_{u}= \squ \vk{b}_{u}-\frac{\vk{x}}{ \squ}}\\
&& \to \ \pk{  \underset{s\in[0,1]}{\underset{ t  \in [0, 1] }\exists }
\begin{array}{ccc}
X_1(t)- \frac{1-\rho}{1+\rho}\mu_1 t >  x_1 \\
X_2(t+s) - \mu (t+s) >x_2
\end{array}
 }
\EQNY
} 
Similarly as before, one can get
\BQNY
Q_{j_1,j_2;u}  \sim   \frac{\widetilde\H(T,S)u^{-1}}{\sqrt{(2\pi)^2\abs{\Sigma_*}}}   e^{-\frac{u}{2} \LT(1+\frac{\frac{(j_2-j_1)T}{4u}}{t^*+\frac{j_1T}{u}+\frac{(j_2-j_1)T}{4u}}\RT) g_L(t^*+\frac{j_1T}{u}+\frac{(j_2-j_1)T}{2u})},
\EQNY
as $u\to\IF,$ where 
\BQNY
\widetilde\H(T,S):=\int_{\R^2}e^{ \vk x^\top \Sigma_{*}^{-1}  \vk{b}_{*} }\pk{
\underset{(t',s')\in\Del_{T,S}}{\underset{ (t,s)\in\Del_{T,S} }\exists }
\begin{array}{ccc}
\frac{1}{2}(X_1(t)+\widetilde X_1(t'))-\frac{\mu_1t^*-1}{4t^*} t-  \frac{\mu_1}{2} t' >  x_1 \\
\frac{1}{2}(X_2(s)+\widetilde X_2(s')) -\frac{\mu_2s^*-1}{4s^*} s- \frac{\mu_2}{2}s' >x_2
\end{array}
 } dx_1 dx_2\in(0,\IF),
\EQNY
with  $(\widetilde X_1, \widetilde X_2)$ an independent copy of $(X_1,X_2)$.
Particularly, letting $j_1=0,j_2=2$ we can show, similarly as \kkk{in (\ref{Lem:Hmumu})}, that
\BQNY
\widetilde\H(T,S) \le \widetilde\H(1,S) (\lfloor T \rfloor )^2 e^{\frac{g_L(t^*)}{8t^*} T}.
\EQNY
Therefore, as $u\to\IF,$
\BQNY
\Pi_{21} (u)
&\lesssim & \frac{\widetilde\H(T,S)u^{-1/2}}{T \sqrt{(2\pi)^2\abs{\Sigma_*}}}  e^{-\frac{u}{2}g_L(t^*)} \LT(\sum^{ N_u^{(1)}} _{ j_1=-N_u^{(1)}} e^{-\frac{b_0}{4}\LT(\frac{j_1 T}{\sqrt u}\RT)^2} \frac{T}{\sqrt u}\RT) \LT(\sum_{j_2=j_1+2}^{N_u^{(1)}} e^{-\frac{g_L(t^*)}{8t^*} (j_2-j_1)T}\RT)\\
&\lesssim& \frac{\widetilde\H(1,S) \lfloor T\rfloor^2  u^{-1/2}}{T \sqrt{(2\pi)^2\abs{\Sigma_*}}}  e^{-\frac{u}{2}g_L(t^*)} \ \int_{\R}e^{-\frac{b_0}{4}   x^2} dx \
 \sum_{k=1}^{\IF} e^{-\frac{g_L(t^*)}{8t^*} kT}. 
\EQNY
Finally, we consider $\Pi_{22}(u)$.
Note \kkk{that}
\BQNY
\Pi_{22}(u) =  \sum^{ N_u^{(1)}} _{ j=- N_u^{(1)}}   p_{j;u}+p_{j+1;u}-\widetilde{p}_{j;u},
\EQNY
where
\BQNY
\widetilde{p}_{j;u}=\pk{\exists_{ ( t,s) \in (t^*+\frac{j T}{u},t^*+\frac{j T}{u})+u^{-1}\Del_{2T,S}} \ X_1(t)>
\sqrt{u}(1+ \mu_1 t ), X_2(s)>
\sqrt{u}(1+ \mu_2 s )}.
\EQNY
Consequently, the claim for $\Pi_{22}(u)$ follows directly by using Lemma \ref{Lem:PL2}. 
\QED

\subsection{Proof of Lemma \ref{Lem:Pr2}} Similarly as in \eqref{eq:DD} we obtain
\BQNY
P_{\theta_0,0}(u)&\le& \pk{
\exists_{(t_1,s_1)\in   \mcE_{11} } \ Z(t_1,s_1) >\sqrt u g(t_1,s_1),\ \
 \exists_{(t_2,s_2)\in  \mcE_{12}}    \    Z(t_2,s_2)>\sqrt u g(t_2,s_2)
}\\
&\le&   \pk{
\exists_{(t_1,s_1)\in   \mcE_{11}  } \   \overline{Z}(t_1,s_1)>\sqrt u\sqrt{g(t_0,s_0)},\ \
 \exists_{(t_2,s_2)\in   \mcE_{12} }    \  \overline{Z}(t_2,s_2)>\sqrt u \sqrt{g(s_0,t_0)}
}\\
&\le&   \pk{
\exists_{(t_1,s_1)\in   \mcE_{11} , (t_2,s_2)\in  \mcE_{12}} \   \overline{Z}(t_1,s_1)+ \overline{Z}(t_2,s_2)> 2\sqrt u \sqrt{g(t_0,s_0)}
},
\EQNY
where, we used the fact that $g(t_0,s_0)=g(s_0,t_0) \le \inf_{(t,s)\in  ( \mcE_{11}\cup \mcE_{12})}g(t,s)$, and
\BQNY
\overline{Z}(t_i,s_i):= \frac{Z(t_i,s_i)}{\sqrt{\mathrm{Var} (Z(t_i,s_i))}}= \frac{Z(t_i,s_i)}{\sqrt{g(t_i,s_i) }}, \ \ \  i=1,2.
\EQNY
Further note that
\BQNY
\E{ Z(t_0,s_0)   Z (s_0,t_0)}&=&\E{(2\mu X_1(t_0) +2(1-2\rho)\mu X_2(s_0)) (2(1-2\rho)\mu X_1(s_0) +2\mu X_2(t_0))}\\
&=& 8(1+2\rho)(1-\rho)\mu.
\EQNY
We obtain
\BQNY
\E{(\overline{Z}(t_0,s_0)+ \overline{Z}(s_0,t_0))^2}&=&2+ 2\E{\overline{Z}(t_0,s_0)  \overline{Z}(s_0,t_0)}\\
&=&2+ 2\frac{\E{ Z(t_0,s_0)   Z (s_0,t_0)}}{g(t_0,s_0)}\\
&=&2+ 2 (1+2\rho)<4.
\EQNY
Thus, for sufficiently small $\theta_0>0$,
\BQNY
\sigma ^2:=\underset{(t_2,s_2)\in     \mcE_{12}}{\underset{(t_1,s_1)\in    \mcE_{11}}{\sup}}\E{(\overline{Z}(t_1,s_1)+ \overline{Z}(t_2,s_2))^2}<4,
\EQNY
where we use continuity of the functions involved. Again, using the Borell-TIS inequality we obtain
\BQNY
\pk{
\exists_{(t_1,s_1)\in     \mcE_{11}, (t_2,s_2)\in     \mcE_{12} } \   \overline{Z}(t_1,s_1)+ \overline{Z}(t_2,s_2)> 2\sqrt u \sqrt{g(t_0,s_0)}}  \le e^{-\frac{(2\sqrt u \sqrt{g(t_0,s_0)}-C_0)^2}{2\sigma ^2}}
\EQNY
holds for all large $u$ such that
\BQNY
2\sqrt u \sqrt{g(t_0,s_0)}>C_0:=\E{  \underset{(t_2,s_2)\in     \mcE_{12}}{\underset{(t_1,s_1)\in     \mcE_{11}}{\sup}} (\overline{Z}(t_1,s_1)+ \overline{Z}(t_2,s_2))   }.
\EQNY
Thus, the claim follows. \QED

\end{appendix}

\section*{Supplemental materials} 

This section 
includes  technical proofs of Lemma A.1 and Lemma A.3.





\subsection*{Proof of Lemma A.1 } \label{a.B}

\COM{
\BEL\label{Lem:Taylor} Assume that $\mu_1<\mu_2$. We have
\begin{itemize}
\item[(i).] If  $-1< \rho  <\hat \rho_1$, then as $(t,s)\to(0,0)$,
\BQNY 
g(t_A+t,s_A+s)=g_A(t_A,s_A)+\frac{a_1}{2}t^2\oo - a_2 ts\oo +\frac{a_3}{2} s^2\oo,
\EQNY
where, with $h(\rho):=\mu_2-2(\mu_1+\mu_2)\rho+3\mu_1\rho^2>0,$
\BQNY
 a_1:=\frac{2\mu_1^3(\mu_2-2\mu_1\rho)}{h(\rho)}>0, \ \
 a_2:=\frac{-2\rho\mu_1^2(\mu_2-2\mu_1\rho)^2}{h(\rho)},  \ \
 a_3: =\frac{2(\mu_2-2\mu_1\rho)^4(1-2\rho)}{h(\rho)}>0.
\EQNY

\item[(ii).] If $ \hat \rho_1  < \rho< \hat \rho_2 $, then
\begin{itemize}
\item{(ii.1),} as $(t,s)\to(0,0)$, with $s<t$ (i.e., $(t^*+t,s^*+s)\in A$),
\BQNY 
g(t^*+t,s^*+s)=g_L(t^*)+b_1( t-s )  \oo + \frac{ c_1 }{2}s^2\oo,
\EQNY
where
\BQNY
b_1&:=&\frac{( \rho-1-2\rho^2)+2\rho(\mu_2-\mu_1\rho)s^* +(1+\rho)\mu_1^2{s^*}^2}{(1-\rho)(1+\rho)^2{s^*}^2}>0,\\
c_1&:=& \frac{2 }{ {s^*}^3}\LT(1+\frac{\rho^2(\rho(1-\rho)-(\mu_2-\mu_1\rho) s^*)^2}{(1-\rho^2)^3}\RT)>0.
\EQNY
\item{(ii.2),} as $(t,s)\to(0,0)$, with $s>t$ (i.e., $(t^*+t,s^*+s)\in B$),
\BQNY 
g(t^*+t,s^*+s)=g_L(t^*)+b_2( s-t )  \oo + \frac{ c_2 }{2}t^2\oo,
\EQNY
where
\BQNY
b_2&:=&\frac{( \rho-1-2\rho^2)+2\rho(\mu_1-\mu_2\rho)t^* +(1+\rho)\mu_2^2{t^*}^2}{(1-\rho)(1+\rho)^2{t^*}^2}>0,  \\
c_2&:=& \frac{2 }{ {t^*}^3}\LT(1+\frac{\rho^2(\rho(1-\rho)-(\mu_1-\mu_2\rho) t^*)^2}{(1-\rho^2)^3}\RT)>0.
\EQNY
\item{(ii.3),} as $(t,s)\to(0,0)$, with $s=t$ (i.e., $(t^*+t,s^*+s)\in L$),
\BQNY 
g(t^*+t, s^*+t)=g_L(t^*)+\frac{b_0}{2}t^2\oo,
\EQNY
where
$b_0:=\frac{4}{(1+\rho){t^*}^3}.
$
\end{itemize}

\item[(iii).]  If $\rho = \hat \rho_1$ (in this case $t_A=s_A=t^*=s^*$), then
  \begin{itemize}
\item{(iii.1),} as $(t,s)\to(0,0)$, with $s<t$,
\BQNY
g(t_A+t,s_A+s)=g_A(t_A,s_A)+\frac{a_1}{2}t^2\oo - a_2 ts\oo +\frac{a_3}{2} s^2\oo,
\EQNY
\item{(iii.2),} as $(t,s)\to(0,0)$, with $s>t$,
\BQNY 
g(t^*+t,s^*+s)=g_L(t^*)+b_2( s-t )  \oo + \frac{ c_2 }{2}t^2\oo.
\EQNY
\item{(iii.3),} as $(t,s)\to(0,0)$, with $s=t$,
\BQNY 
g(t^*+t, s^*+t)=g_L(t^*)+\frac{b_0}{2}t^2\oo.
\EQNY

\end{itemize}


\end{itemize}

\EEL

} 

(i). Recall that for any $-1< \rho< \hat\rho_1,$ the global minimizer is given by
\BQNY
(t_A,s_A)=   \LT( \frac{ 1-2\rho}{  \mu_1}, \frac{1}{ \mu_2-2\mu_1 \rho } \RT)\in A.
\EQNY
In a small neighbourhood of  this point, we have
$$
g(t,s)=g_A(t,s)=\frac{(1+\mu_2 s)^2}{s} +\frac{((1+\mu_1 t)-\rho(1+\mu_2 s))^2}{t-\rho^2 s}.
$$
Next we calculate the first few partial derivatives.  We have
\BQNY 
&&\frac{\partial g(t,s)}{\partial t} =\frac{ (\mu_1 t+1-\rho-\rho\mu_2 s)(\mu_1 t -(2\mu_1\rho^2-\rho \mu_2)s+\rho-1)}{(t-\rho^2s)^2}\nonumber\\
&&\frac{\partial g(t,s)}{\partial s} = -\frac{1}{s^2}+\mu_2^2-2\mu_2\rho
\frac{ \mu_1 t+1-\rho-\rho\mu_2 s}{t-\rho^2s}+\rho^2 \LT(\frac{ \mu_1 t+1-\rho-\rho\mu_2 s}{t-\rho^2s}\RT)^2.
\EQNY
Since
\BQNY
\mu_1 t_A -(2\mu_1\rho^2-\rho \mu_2)s_A+\rho-1=0,\ \ \ \frac{ \mu_1 t_A+1-\rho-\rho\mu_2 s_A}{t_A-\rho^2s_A}=2\mu_1,
\EQNY
one can obtain that
\BQNY
\frac{\partial g(t,s)}{\partial t}  \mid _{(t_A,s_A)}  =\frac{\partial g(t,s)}{\partial s}\mid _{(t_A,s_A)}  =0.
\EQNY
Further calculations show that
\BQN\label{eq:a1}
&&\frac{\partial^2 g(t,s)}{\partial t^2}  \mid _{(t_A,s_A)}=  \frac{2\mu_1^2}{t-\rho^2 s} \mid _{(t_A,s_A)}  =\frac{2\mu_1^3(\mu_2-2\mu_1\rho)}{\mu_2-2(\mu_1+\mu_2)\rho+3\mu_1\rho^2}=a_1>0, \\
&& \frac{\partial^2 g(t,s)}{\partial t \partial s}  \mid _{(t_A,s_A)} =  \frac{2\mu_1 \rho (\mu_2-2\mu_1\rho)}{t-\rho^2 s} \mid _{(t_A,s_A)}  =\frac{2\rho\mu_1^2(\mu_2-2\mu_1\rho)^2}{\mu_2-2(\mu_1+\mu_2)\rho+3\mu_1\rho^2}=-a_2,  \nonumber\\
&&\frac{\partial^2 g(t,s)}{\partial s^2}  \mid _{(t_A,s_A)}  =\frac{2(\mu_2-2\mu_1\rho)^4(1-2\rho)}{\mu_2-2(\mu_1+\mu_2)\rho+3\mu_1\rho^2}= a_3>0.\nonumber
\EQN
Consequently, the claim in (i) is established.

(ii). For $\hat\rho_1<\rho<\hat\rho_2$,
 the  global minimizer is given by
\BQNY
(t^*,s^*) \in L.
\EQNY
Three different types of expansion are available for $g(t^*+t,s^*+s)$, according to
\BQNY
(ii.1). (t^*+t,s^*+s)\in A,\ \ \ (ii.2). (t^*+t,s^*+s)\in B,\ \ \ (ii.3). (t^*+t,s^*+s)\in  L.
\EQNY
\underline{(ii.1) $(t^*+t,s^*+s)\in A$.} Consider the  two representations of $g_3(t,s)$ given in (12) and (13). 
\COM{
\BQN
g_A(t,s)&=&\frac{(1+\mu_1 t)^2s -2\rho s(1+\mu_2 s)(1+\mu_1 t)+(1+\mu_2 s)^2 t}{ts-\rho^2s^2}
&=&\frac{(1+\mu_2 s)^2}{s} +\frac{((1+\mu_1 t)-\rho(1+\mu_2 s))^2}{t-\rho^2 s}.
\EQN
}
Using the representation (12) 
we can show that
\BQNY
\frac{\partial g_A(t,s)}{\partial t}  \mid _{(t^*,s^*)}  =-\frac{\partial g_A(t,s)}{\partial s}\mid _{(t^*,s^*)}  =\frac{( \rho-1-2\rho^2)+2\rho(\mu_2-\mu_1\rho)s^* +(1+\rho)\mu_1^2{s^*}^2}{(1-\rho)(1+\rho)^2{s^*}^2}=b_1.
\EQNY
On the other hand, using representation in (13) 
we can show that $b_1>0,$ i.e.,
\BQNY
b_1= \frac{\partial g_A(t,s)}{\partial t}  \mid _{(t^*,s^*)} =\frac{((\mu_1-\rho\mu_2)s^*+1-\rho)((\mu_1-2\mu_1\rho^2+\rho\mu_2)s^*+\rho-1)}{{s^*}^2 (1-\rho^2)^2}>0.
\EQNY
In fact,  for $\rho\in(0,\mu_1/\mu_2]$ we have that
\BQNY
(\mu_1-\rho\mu_2)s^*+1-\rho>0,
\EQNY
and for $\rho\in(\mu_1/\mu_2, \hat\rho_2)$ we have by (d).(i) of  Lemma 7 in \cite{DJT19a}  (note $s^*=t^*$) that the above still holds. Furthermore, we have by (b).(ii) of Lemma 9 in \cite{DJT19a} that, for all $\rho\in(\hat\rho_1, 1)$
\BQNY
(\mu_1-2\mu_1\rho^2+\rho\mu_2)s^*+\rho-1>0.
\EQNY
Moreover, by using representation (13) 
one can show  that
\BQNY
\frac{\partial^2 g_A(t,s)}{\partial s^2}  \mid _{(t^*,s^*)}  =\frac{2 }{ {s^*}^3}\LT(1+\frac{\rho^2(\rho(1-\rho)-(\mu_2-\mu_1\rho) s^*)^2}{(1-\rho^2)^3}\RT)= c_1>0.
\EQNY
Then, for  $(t^*+t,s^*+s)\in A$ (where $t>s$), we have by Taylor expansion
\BQNY 
g(t^*+t,s^*+s)&=&g_A(t^*,s^*)+b_1( t-s )  \oo + \frac{1 }{2}\frac{\partial^2 g_A(t,s)}{\partial t^2} \mid _{(t^*,s^*)}   t^2\oo \nonumber\\
&&\ \ \ +  \frac{\partial^2 g_A(t,s)}{\partial t \partial s} \mid _{(t^*,s^*)}  ts\oo +\frac{ c_1 }{2}s^2\oo \nonumber\\
&=&g_L(t^*)+b_1( t-s )  \oo + \frac{ c_1 }{2}s^2\oo,\ \  (t,s)\to(0,0)
\EQNY
as required.

\underline{(ii.2) $(t^*+t,s^*+s)\in B$.}
Similarly as in (ii.1), we have
\BQNY
\frac{\partial g_B(t,s)}{\partial s}  \mid _{(t^*,s^*)}  =-\frac{\partial g_B(t,s)}{\partial t}\mid _{(t^*,s^*)}  =\frac{( \rho-1-2\rho^2)+2\rho(\mu_1-\mu_2\rho)t^* +(1+\rho)\mu_2^2{t^*}^2}{(1-\rho)(1+\rho)^2{t^*}^2}=b_2>0,
\EQNY
where $b_2>0$ follows similarly as the positiveness of $b_1$, by using Lemma 7 in  \cite{DJT19a}.
Moreover,  one can show  that
\BQNY
\frac{\partial^2 g_B(t,s)}{\partial t^2}  \mid _{(t^*,s^*)}  =\frac{2 }{ {t^*}^3}\LT(1+\frac{\rho^2(\rho(1-\rho)-(\mu_1-\mu_2\rho) t^*)^2}{(1-\rho^2)^3}\RT)= c_2>0.
\EQNY
Consequently, by Taylor expansion the claim of (ii.2) is established.
\COM{
\BQN\label{eq:Pg2}
g_B(t^*+t,s^*+s)=g_B(t^*,s^*)+b_2( s-t )  \oo + \frac{ c_2 }{2}t^2\oo, \ \  (t,s)\to(0,0).
\EQN
}

\underline{(ii.3) $(t^*+t,s^*+s)\in L$.} The claim follows by considering the univariate function $g_L(s)$ defined in (14). 
\COM{In this case, consider
\BQNY
f(t)=g(t,t)=\frac{(1+\mu_1 t)^2-2\rho(1+\mu_1t)(1+\mu_2 t)+(1+\mu_2 t)^2 }{(1-\rho^2)t}.
\EQNY
Further, for any $(t^*+t,s^*+t)\in  L,$
\BQN\label{eq:Pg3}
g(t^*+t,s^*+t)=g(t^*,s^*)+\frac{b_0}{2}t^2\oo, \ \ \  t\to 0,
\EQN
where
\BQNY
b_0=\frac{4}{(1+\rho){t^*}^3}.
\EQNY
}

(iii).  
Note that in this case $b_1=0$ and $b_2>0$. The claims follows by combining the above discussions in (i)-(ii). This completes the proof.

\subsection*{Proof of Lemma A.3}

\COM{
\BEL\label{Lem:fPH}
We have, for any $T,S>0$
\BQNY 
\lim_{u\to\IF}\int_{\R^{2}}  f_{j,l;u}(\vk x)P_{j,l;u}(\vk x)\,d\vk{x}=\H (\mu_1;T)\H(\mu_2-2\mu_1\rho;S)
\EQNY
holds uniformly for  $-N_u^{(1)}\le j\le N_u^{(1)}, -N_u^{(2)}\le l\le N_u^{(2)}.$
\EEL
}

 The proof consists of two steps. In Step I, we derive the limit, as $u\to\IF,$ of the integrand $f_{j,l;u}(\vk x)P_{j,l;u}(\vk x)$ for any fixed $\vk x$. In Step II, we look for uniform integrable upper bound of $f_{j,l;u}(\vk x)P_{j,l;u}(\vk x)$, by which we show that the limit can pass into the integral.  For simplicity, in the following when we write $j,l$, we mean   $-N_u^{(1)}\le j\le N_u^{(1)}, -N_u^{(2)}\le l\le N_u^{(2)}.$

\underline{Step I.} Recall
\BQNY
 P_{j,l;u}(\vk x) = \pk{\sup_{t\in[0,T]}\  X_1(t ) -\mu_1 t>     x_1 }\  \pk{\sup_{s\in[0,S]} \ Y_{l;u}(s)> x_2  \Big{|} \vk Z_{j,l;u}= \squ \vk{b}_{j,l;u}-\frac{\vk{x}}{ \squ}},
\EQNY
 where
 \BQNY
 Y_{l;u}(s)= \sqrt u \LT(X_2(s_A+\frac{lS}{u}+\frac{s}{u})- X_2(s_A+\frac{lS}{u})\RT) -  \mu_2 s.
 \EQNY
It follows that $( Y_{l;u}(s), \vk Z_{j,l;u}^\top)$ is normally distributed with mean value and covariance matrix given, respectively, by
\BQNY
\vk m_{j,l;u}(s)=(-\mu_2 s, 0, 0), \ \ \
\widehat\Sigma_{j,l;u}(s)=
\left(
   \begin{array}{ccr}
   s&  \frac{\rho s}{\sqrt u} & 0\\
 \frac{\rho s}{\sqrt u}   &  t_A+\frac{jT}{u}  & \rho\  (s_A+\frac{lS}{u})\\
  0&    \rho\  (s_A+\frac{lS}{u}) & s_A+\frac{lS}{u} \\
   \end{array}
 \right).
\EQNY
 Therefore, denoting for any $\vk x\in \R^2$
 \BQN\label{eq:WW}
 W_{j,l;u}(s)= Y_{l;u}(s)  \mid \vk Z_{j,l;u}= \squ \vk{b}_{j,l;u}-\frac{\vk{x}}{ \squ},
 \EQN
we can show that its mean   is given by
\BQN\label{eq:mean}
\E{ W_{j,l;u}(s)}&=& -\mu_2 s + (\frac{\rho s}{\sqrt u}, 0)\ (\Sigma_{i,j;u})^{-1} \ \LT( \squ \vk{b}_{j,l;u}-\frac{\vk{x}}{ \squ}\RT)\nonumber\\
 &=&-\mu_2 s + \frac{(a_{j;u}-\rho b_{l;u})-\frac{x_1-\rho x_2}{u}}{(t_0+\frac{jT}{u}) -\rho^2(s_0+\frac{lS}{u})}\ \rho s,
\EQN
and its variance  is given by
\BQN\label{eq:var}
\mathrm{Var}(W_{j,l;u}(s)) &=&  s - (\frac{\rho s}{\sqrt u}, 0)\  \left(
  \Sigma_{i,j;u}^{-1}
 \right)^{-1}\  (\frac{\rho s}{\sqrt u}, 0)^\top\nonumber\\
 &=& s \LT( 1- \frac{ \rho^2 }{(t_0+\frac{jT}{u}) -\rho^2(s_0+\frac{lS}{u})}\ \frac{ s}{u} \RT).
\EQN
Similarly, we can derive that, for any $s_1,s_2\in[0,S],$
\BQN\label{eq:Cov}
\mathrm{Var}(  W_{j,l;u}(s_1) -  W_{j,l;u}(s_2)) =\abs{s_1-s_2} \LT( 1- \frac{ \rho^2 }{(t_0+\frac{jT}{u}) -\rho^2(s_0+\frac{lS}{u})}\ \frac{ \abs{s_1-s_2}}{u} \RT).
\EQN
By \eqref{eq:mean} and \eqref{eq:Cov} we have, for any $s, s_1,s_2\in[0,S]$,
\BQN
&&\E{ W_{j,l;u}(s)} \to -( \mu_2-2\mu_1\rho ) s <0, \label{eq:meanW}\\
&&\mathrm{Var}(  W_{j,l;u}(s_1) -  W_{j,l;u}(s_2)) \to \abs{s_1-s_2}=\mathrm{Var}(  B_1(s_1) -  B_2(s_2)) \nonumber
\EQN
as $u\to\IF$,
where the convergence is uniform with respect to $j,l$. 
Thus, by Lemma 4.2 in \cite{ZhouXiao17} we conclude that the finite dimensional distributions of $W_{j,l;u}(s), s\in[0,S]$ converge to the finite dimensional distributions of $B_1(s) -( \mu_2-2\mu_1\rho ) s, s\in[0,T]$ as $u\to\IF$ uniformly with respect to $j, l$. 
Further, note that from \eqref{eq:Cov} we have  
\BQNY
\mathrm{Var}(  W_{j,l;u}(s_1) -  W_{j,l;u}(s_2)) \le   \abs{s_1-s_2}
\EQNY
holds, for all $j, l$, 
when $u$ is large. This guarantees the uniform tightness of  $W_{j,l;u}(s), s\in[0,S]$ (see, e.g., Proposition 9.7 in \cite{Pit15}), and thus we conclude that the stochastic processes
$W_{j,l;u}(s), s\in[0,S]$ converge weakly to $B_1(s) -( \mu_2-2\mu_1\rho )s, s\in[0,S]$ as $u\to\IF$ uniformly with respect to $ j, l$. Therefore,  we have from the continuous mapping theorem that for any $\vk x\in \R^2$
\BQNY
\pk{\sup_{s\in[0,S]} \ Y_{l;u}(s)> x_2  \Big{|} \vk Z_{j,l;u}= \squ \vk{b}_{j,l;u}-\frac{\vk{x}}{ \squ}} \to 
\pk{\sup_{s\in[0,S]} B_1(s) -( \mu_2-2\mu_1\rho ) s>x_2}
\EQNY
 as $u\to\IF$, uniformly with respect to $j,l$. 
Further, it follows that
\BQNY
\Sigma_{j,l;u}^{-1}  \vk{b}_{j,l;u} \to (2\mu_1, 2( \mu_2-2\mu_1\rho ))^\top, \ \ u\to\IF
\EQNY
holds uniformly with respect to $ j, l$. Thus,
\BQNY
f_{j,l;u}(\vk x) \to  \exp\LT( 2\mu_1 x_1+ 2( \mu_2-2\mu_1\rho ) x_2 \RT), \ \ \ \ u\to\IF
\EQNY
holds uniformly with respect to $j, l$.

\underline{Step II.} In order to pass the limit into the integral, it is sufficient to find an  integrable upper bound $h(\vk x)$ such that
\BQNY
f_{j,l;u}(\vk x) P_{j,l;u}(\vk x)\le h(\vk x)
\EQNY
holds for all large $u$, uniformly with respect to $j, l$. The four quadrants will be considered separately.

\underline{(i).  $x_1<0, x_2<0$.} In this case, an upper bound for $P_{j,l;u}(\vk x)$ is chosen to be 1, and for some small $\vn>0$
\BQNY
f_{j,l;u}(\vk x) \le \exp\LT( 2\mu_1 (1-\vn) x_1+ 2(\mu_2-2\mu_1\rho)  (1-\vn) x_2 \RT).
\EQNY
Thus, we choose
\BQNY
h(\vk x)=\exp\LT( 2\mu_1 (1-\vn) x_1+  2(\mu_2-2\mu_1\rho)   (1-\vn) x_2 \RT).
\EQNY

\underline{(ii).  $x_1>0, x_2<0$.}  In this case,
\BQNY
 P_{j,l;u}(\vk x) \le \pk{\sup_{t\in[0,T]}\  X_1(t ) >     x_1 }=\pk{\abs{X_1(T)}>x_1}\le \frac{2}{\sqrt{2\pi T} x_1}e^{-\frac{x_1^2}{2T}}
\EQNY
and for some small $\vn>0$
\BQNY
f_{j,l;u}(\vk x) \le \exp\LT( 2\mu_1 (1+\vn) x_1+  2(\mu_2-2\mu_1\rho)  (1-\vn) x_2 \RT).
\EQNY
 Thus, for any $M>0$ we choose (with $I_{(\cdot)}$ denoting the indicator function)
 \BQNY
h(\vk x)=\exp\LT( 2\mu_1 (1+\vn) x_1+ 2(\mu_2-2\mu_1\rho)  (1-\vn) x_2 \RT) \LT(I_{(x_1<M)}+I_{(x_1\ge M)}\frac{2}{\sqrt{2\pi T} x_1}e^{-\frac{x_1^2}{2T}} \RT).
\EQNY

 \underline{(iii).  $x_1<0, x_2>0$.}  In this case,
 for some small $\vn>0$
\BQNY
f_{j,l;u}(\vk x) \le \exp\LT(  2\mu_1 (1-\vn)  x_1 + 2(\mu_2-2\mu_1\rho)   (1+\vn) x_2 \RT),
\EQNY
and
\BQNY
 P_{j,l;u}(\vk x) \le \pk{\sup_{s\in[0,S]} \ W_{j,l;u}(s)> x_2  }
\EQNY

\COM{
and

 Thus, for any $M>0$ we choose
 \BQNY
g(\vk x)=\exp\LT( 2\mu (1+\vn) x_1+ 2(1-2\rho)\mu  (1-\vn) x_2 \RT) \LT(I_{(x_1<M)}+I_{(x_1\ge M)}\frac{2}{\sqrt{2\pi T} x_1}e^{-\frac{x_1^2}{2T}} \RT).
\EQNY

 Next, for any sufficiently large constant $M>0$ we write
 \BQNY
&&\int_{\R^{2}}  f_{j,l;u}(\vk x)P_{j,l;u}(\vk x)\,d\vk{x}\\&&= \int_{\abs{x_1}\le M, \abs{x_2}\le M } + \int_{\abs{x_1}\le M, \abs{x_2}> M } + \int_{\abs{x_1}> M, \abs{x_2}\le M } +\int_{\abs{x_1}> M, \abs{x_2}> M }  f_{j,l;u}(\vk x)P_{j,l;u}(\vk x)\,d\vk{x} \\
&&= I_1(M,j,l;u)+I_2(M,j,l;u)+I_3(M,j,l;u)+I_4(M,j,l;u)
\EQNY

Thus, by the bounded convergence theorem   we conclude that, as $u\to\IF,$
\BQNY
 I_1(M,j,l;u) &\to& \int_{\abs{x_1}\le M }  e^{2\mu  x_1 }  \pk{\sup_{t\in[0,T]}\  X_1(t ) -\mu_1 t>     x_1 } dx_1\\
 && \times \int_{\abs{x_2}\le M }  e^{2(1-2\rho )\mu  x_2 } \pk{\sup_{s\in[0,S]} B_1(s) -(1-2\rho )\mu_2 s>x_2}  dx_2
\EQNY
Now we look at $I_k(M,j,l;u), k=1,2,3,$ respectively. Note that
\BQNY
 I_2(M,j,l;u) &= & \int_{\abs{x_1}\le M, x_2<-M }   + \int_{\abs{x_1}\le M, x_2 > M }  f_{j,l;u}(\vk x)  P_{j,l;u}(\vk x)  dx_1 dx_2 =:  I_{2,1}(M,j,l;u) +  I_{2,2}(M,j,l;u)
\EQNY
We have, for all $\abs{x_1}\le M, x_2<-M,$ and all large $u$
\BQNY
 f_{j,l;u}(\vk x)  P_{j,l;u}(\vk x) \le   f_{j,l;u}(\vk x) \pk{\sup_{t\in[0,T]}\  X_1(t ) -\mu_1 t>     x_1 }  \le
 C(\vn, M)\  e^{2(1-2\rho )\mu (1-\vn)  x_2}
\EQNY
holds, with $\vn\in(0,1)$ and $C(\vn, M)$ some positive constant  depending on $\vn,M$. Then, by the dominated  convergence theorem
\BQNY
I_{2,1}(M,j,l;u) &\le&  \int_{\abs{x_1}\le M, x_2 <- M }    f_{j,l;u}(\vk x)   \pk{\sup_{t\in[0,T]}\  X_1(t ) -\mu_1 t>     x_1 } dx_1\\
&\to&  \int_{\abs{x_1}\le M }  e^{2\mu  x_1 }  \pk{\sup_{t\in[0,T]}\  X_1(t ) -\mu_1 t>     x_1 } dx_1  \int_{ x_2 <- M }  e^{2(1-2\rho )\mu  x_2 }   dx_2.
\EQNY
Further, we have (recall that $W_{j,l;u}(s)$ depends also on $x_1,x_2$; see \eqref{eq:WW})
\BQN \label{eq:I22}
 I_{2,2}(M,j,l;u) &= & \int_{\abs{x_1}\le M,  x_2 > M}  
  f_{j,l;u}(\vk x)   \pk{\sup_{t\in[0,T]}\  X_1(t ) -\mu_1 t>     x_1 } \nonumber \\
  &&\ \qquad\qquad \times  \pk{\sup_{s\in[0,S]} \ W_{j,l;u}(s)> x_2  } dx_1 dx_2.
\EQN
}

Next we derive an upper bound for $\pk{\sup_{s\in[0,S]} \ W_{j,l;u}(s)> x_2  }$.
It follows from \eqref{eq:mean} and \eqref{eq:meanW} that 
\BQNY
\sup_{s\in{[0,S]}}\E{W_{j,l;u}(s)} \le \frac{  (\rho^2x_2-\abs{\rho} x_1)   /{u}}{(t_A+\frac{jT}{u}) -\rho^2(s_A+\frac{lS}{u})}\   S \le  \vn_0 x_2+c_0\abs{x_1}/u
\EQNY
holds for all $u$ large and uniformly in $ j, l$, with some small constant $\vn_0>0$ and $c_0>0$. Then
\BQNY
 \pk{\sup_{s\in[0,S]} \ W_{j,l;u}(s)> x_2  } \le  \pk{\sup_{s\in[0,S]} \ W_{j,l;u}(s)-\E{W_{j,l;u}(s)}> (1-\vn_0)x_2 - c_0\abs{x_1} /u}.
\EQNY
We consider the following two subsets of $E:=\{(x_1,x_2): x_1<0, x_2>0\}$:
$$
E_1:=\{(x_1,x_2): (1-\vn_0)x_2 - c_0\abs{x_1} /u > \vn_1 x_2\}\cap E,\ \  E_2:=\{(x_1,x_2): (1-\vn_0)x_2 - c_0\abs{x_1} /u \le  \vn_1 x_2\}\cap E,
$$
with $\vn_1>0$ some small constant.

Below, we derive upper bounds $h(\vk x)$ on $E_1, E_2 $,
respectively.
Note,  from \eqref{eq:Cov},
for all large $u$ and  all $ j, l$,
\BQNY
 \mathrm{Var}(  W_{j,l;u}(s_1) -  W_{j,l;u}(s_2)) \le \mathrm{Var}(  B_1(s_1) -  B_2(s_2))
\EQNY
for any $s_1,s_2\in[0,S]$. Hence, by the Sudakov-Fernique inequality (see, e.g., \cite{AdlerTaylor})
\BQNY
 \E{\sup_{s\in [0,S]} (W_{j,l;u}(s) - \E{W_{j,l;u}(s)})} \le  \E{\sup_{s\in [0,S]} B_1(s)}:=U_0.
\EQNY
Moreover, it follows from \eqref{eq:var} that, for sufficiently large $u$ and  all $j, l$,
\BQNY
\sup_{s\in [0,S]} \mathrm{Var}(  W_{j,l;u}(s))\le S.
\EQNY
Thus, we have from the Borell-TIS inequality (Theorem 2.1.1 in \cite{AdlerTaylor})
that, on $E_1,$ for all $u$ large and $x_2>{U_0}/{\vn_1}$ 
\BQNY
 \pk{\sup_{s\in[0,S]} \ W_{j,l;u}(s)> x_2  }&\le& \pk{\sup_{s\in[0,S]} \ W_{j,l;u}(s)-\E{W_{j,l;u}(s)}>  \vn_1x_2  }\\
 &\le & \exp\LT(-\frac{(\vn_1 x_2 -U_0)^2}{2S}\RT).
\EQNY
Therefore,  on $E_1,$ we can choose
 \BQNY
h(\vk x)&=&\exp\LT( 2\mu (1-\vn) x_1+ 2(1-2\rho)\mu  (1+\vn) x_2 \RT)\nonumber \\
&&\times \LT(I_{(x_2\le  {U_0}/\vn_1)}+I_{(x_2>{U_0}/\vn_1)}  \exp\LT(-\frac{( \vn_1 x_2 -U_0)^2}{2S}\RT) \RT).\label{eq:gT}
\EQNY
On the other hand,  we have,  for some  $M> \frac{2(\mu_2-2\mu_1\rho)   (1+\vn) }{(1-\vn_0-\vn_1)\mu_1 (1-\vn)} $,
\BQNY
E_2  \subseteq E_3:=\{  (x_1,x_2):   (1-\vn_0-\vn_1)x_2 \le c_0\abs{x_1}/u \le  \abs{x_1}/M\}\cap E,
\EQNY
holds for all large $u$. Thus, on $E_2$ we can choose 
\BQNY
h(\vk x) = \exp\LT( \mu_1 (1-\vn)  x_1 - (M(1-\vn_0-\vn_1)\mu_1 (1-\vn)- 2(\mu_2-2\mu_1\rho)   (1+\vn)) x_2 \RT).
\EQNY

 \underline{(iv).  $x_1>0, x_2>0$.}  In this case,
 for some small $\vn>0$
\BQNY
f_{j,l;u}(\vk x) \le \exp\LT( 2\mu (1+\vn) x_1+ 2(1-2\rho)\mu  (1+\vn) x_2 \RT),
\EQNY
and
\BQNY
 P_{j,l;u}(\vk x) \le  \pk{\sup_{s\in[0,S]} \ W_{j,l;u}(s)> x_2  } \pk{\abs{X_1(T)}>x_1}.
\EQNY
\COM{ 
Let $M>\frac{U_0}{(1-\vn_0) \vn_1}$ be some constant, with $\vn_1$  and $\vn_0$ the same as above. We consider the following three sub-scenarios:
(S.1). $M>x_1>0, x_2>0$; (S.2). $x_1>M, M>x_2>0$; (S.3). $x_1>M, x_2>M$.

\underline{(S.1). $M>x_1>0, x_2>0$}. Similarly to (iii) we choose
 \BQNY
g(\vk x)&=&\exp\LT( 2\mu (1+\vn) x_1+ 2(1-2\rho)\mu  (1+\vn) x_2 \RT)\nonumber \\
&&\times \LT(I_{(x_2\le \frac{U_0}{1-\vn_0})}+I_{(x_2>\frac{U_0}{1-\vn_0})}  \exp\LT(-\frac{((1-\vn_0) x_2 -U_0)^2}{2T}\RT) \RT).\label{eq:gT}
\EQNY

\underline{(S.2). $x_1>M, M>x_2>0$}. Similarly to (ii) we choose
 \BQNY
g(\vk x)=\exp\LT( 2\mu (1+\vn) x_1+ 2(1-2\rho)\mu  (1+\vn) x_2 \RT) \frac{2}{\sqrt{2\pi T} x_1}e^{-\frac{x_1^2}{2T}}.
\EQNY

\underline{(S.3). $x_1>M, x_2>M$}.
First recall that  
\BQNY
 \pk{\abs{X_1(T)}>x_1} \le \frac{2}{\sqrt{2\pi T} x_1}e^{-\frac{x_1^2}{2T}}.
\EQNY
Next, we derive an upper bound for
$$
\pk{\sup_{s\in[0,S]} \ W_{j,l;u}(s)> x_2  }
$$
Similarly as before we have
\BQNY
\sup_{s\in{[0,T]}}\E{W_{j,l;u}(s)} \le \frac{ \frac{\rho^2 x_2 -\rho  x_1  }{u}}{(t_0+\frac{jT}{u}) -\rho^2(s_0+\frac{lS}{u})} T  \le \vn_0 x_2  -\rho a_0  x_1 /u
\EQNY
for all $u$ large and uniformly in $-N_u\le j,l\le N_u$, with some $a_0>0$. 
Then, with $c_0:=  - a_0 \rho >0,$
\BQNY
 \pk{\sup_{s\in[0,S]} \ W_{j,l;u}(s)> x_2  }&\le& \pk{\sup_{s\in[0,S]} \ W_{j,l;u}(s)-\E{W_{j,l;u}(s)}> (1-\vn_0)x_2 - c_0x_1 /u}.
\EQNY
Consider the following two subsets of $E:=\{(x_1,x_2): x_1>M, x_2>M\}$:
$$A:=\{(x_1,x_2): (1-\vn_0)x_2 - c_0x_1 /u > \vn_1 x_2\}\cap E,\ \ \ B:=\{(x_1,x_2): (1-\vn_0)x_2 - c_0x_1 /u \le  \vn_1 x_2\}\cap E.$$
Note that on the subset $A$ we have by the Borell-TIS inequality
\BQNY
 \pk{\sup_{s\in[0,S]} \ W_{j,l;u}(s)> x_2  }&\le& \pk{\sup_{s\in[0,S]} \ W_{j,l;u}(s)-\E{W_{j,l;u}(s)}>  \vn_1 x_2}\\
 &\le& \exp\LT(-\frac{  (\vn_1 x_2 -U_0)^2}{2T}\RT).
\EQNY
On the other hand   we have for large $u$
\BQNY
B  \subseteq D:=\{  (x_1,x_2):   (1-\vn_0-\vn_1)x_2 \le c_0x_1/u \le x_1\}\cap E.
\EQNY
Thus, we obtain
\BQNY
 P_{j,l;u}(\vk x) f_{j,l;u}(\vk x) & \le& \exp\LT(2\mu (1+\vn)  x_1 +  2(1-2\rho )\mu  (1+\vn) x_2  \RT) \frac{2}{\sqrt{2\pi T} x_1} \\
 &&\times \LT(\exp\LT(-\frac{x_1^2}{2T}-\frac{  (\vn_1 x_2 -U_0)^2}{2T}\RT) I_{A}+\exp\LT(-\frac{x_1^2}{4T}- \frac{(1-\vn_0-\vn_1)^2x_2^2}{4T}\RT)I_{D}  \RT).
\EQNY
}
By similar arguments as (iii), we can  choose
\BQNY
h(\vk x)&=&\exp\LT(2\mu (1+\vn)  x_1 +  2(1-2\rho )\mu  (1+\vn) x_2  \RT) \frac{2}{\sqrt{2\pi T} x_1} \\
 &&\times \LT(\exp\LT(-\frac{x_1^2}{2T}-\frac{  (\vn_1 x_2 -U_0)^2}{2S} I_{(x_2>U_0/\vn_1)}\RT) +\exp\LT(-\frac{x_1^2}{4T}- \frac{(1-\vn_0-\vn_1)^2x_2^2}{4T}\RT)   \RT).
\EQNY
This competes the proof.

\bigskip
{\bf Acknowledgement}:
TR \& KD were partially supported by NCN Grant No  2018/31/B/ST1/00370 (2019-2022).

\bibliographystyle{plain}

 \bibliography{vectProcEKEE}

\def\polhk#1{\setbox0=\hbox{#1}{\ooalign{\hidewidth
  \lower1.5ex\hbox{`}\hidewidth\crcr\unhbox0}}}
  \def\lfhook#1{\setbox0=\hbox{#1}{\ooalign{\hidewidth
  \lower1.5ex\hbox{'}\hidewidth\crcr\unhbox0}}}
  \def\polhk#1{\setbox0=\hbox{#1}{\ooalign{\hidewidth
  \lower1.5ex\hbox{`}\hidewidth\crcr\unhbox0}}}
\begin{thebibliography}{10}

\bibitem{AdlerTaylor}
R.J. Adler and J.E. Taylor.
\newblock {\em Random fields and geometry}.
\newblock Springer Monographs in Mathematics. Springer, New York, 2007.

\bibitem{AzW09}
J.~Aza{\"\i}s and M.~Wschebor.
\newblock {\em Level sets and extrema of random processes and fields}.
\newblock John Wiley \& Sons, 2009.

\bibitem{DHJT18}
K.~D\c{e}bicki, E.~Hashorva, L.~Ji, and T.~Rolski.
\newblock Extremal behavior of hitting a cone by correlated {B}rownian motion
  with drift.
\newblock {\em Stochastic Processes and their Applications},
  128(12):4171--4206, 2018.

\bibitem{KEP2015}
K.~D\c{e}bicki, E.~Hashorva, and P.~Liu.
\newblock Ruin probabilities and passage times of $\gamma$-reflected {G}aussian
  process with stationary increments.
\newblock {\em ESAIM: Probability and Statistics}, 21:495--535, 2017.

\bibitem{DHM19}
K.~{D{\c{e}}bicki}, E.~{Hashorva}, and Z.~{Michna}.
\newblock {Simultaneous Ruin Probability for Two-Dimensional Brownian Risk
  Model}.
\newblock {\em Accepted for publication in J. Appl. Probab.}, 2019.

\bibitem{DHW19}
K.~{D{\c{e}}bicki}, E.~{Hashorva}, and L.~{Wang}.
\newblock {Extremes of vector-valued Gaussian processes}.
\newblock {\em arXiv e-prints: arXiv:1911.06350}, 2019.

\bibitem{DJT19a}
K.~D\c{e}bicki, L.~Ji, and T.~Rolski.
\newblock Logarithmic asymptotics for probability of component-wise ruin in
  two-dimensional {B}rownian model.
\newblock {\em Risks}, 7(83), 2019.

\bibitem{RolskiSPA}
K.~D{\polhk{e}}bicki, K.~M. Kosi{\'n}ski, M.~Mandjes, and T.~Rolski.
\newblock Extremes of multidimensional {G}aussian processes.
\newblock {\em Stochastic Process. Appl.}, 120(12):2289--2301, 2010.

\bibitem{HJ14c}
E.~Hashorva and L.~Ji.
\newblock Extremes and first passage times of correlated fractional {B}rownian
  motions.
\newblock {\em Stochastic Models}, 30(3):272--299, 2014.

\bibitem{HKR98}
H.~He, W.~P. Keirstead, and J.~Rebholz.
\newblock Double lookbacks.
\newblock {\em Mathematical Finance}, 8(3):201--228, 1998.

\bibitem{Honnappa}
H.~Honnappa, P.~Jaiswal, and R.~Pasupathy.
\newblock Large deviations of gaussian extremes on convex sets.
\newblock {\em Manuscript.
  https://web.ics.purdue.edu/~pasupath/PAPERS/ldextremes.pdf}, 2020.

\bibitem{Ji18}
L.~Ji.
\newblock On the cumulative {P}arisian ruin of multi-dimensional {B}rownian
  motion models.
\newblock {\em Preprint, https://arxiv.org/pdf/1811.10110.pdf}, 2019.

\bibitem{KZ16}
S.~Kou and H.~Zhong.
\newblock First-passage times of two-dimensional {B}rownian motion.
\newblock {\em Adv. Appl. Prob.}, 48:1045--1060, 2016.

\bibitem{LM07}
P.~Lieshout and M.~Mandjes.
\newblock Tandem {B}rownian queues.
\newblock {\em Math. Methods Oper. Res.}, 66:275--298, 2007.

\bibitem{Mandjes}
M.~Mandjes.
\newblock {\em Large {D}eviations for {G}aussian {Q}ueues: {M}odelling
  {C}ommunication {N}etworks.}
\newblock Wiley, Chichester, 2007.

\bibitem{MO67}
A.~W. Marshall and I.~Olkin.
\newblock A multivariate exponential distribution.
\newblock {\em J.Amer. Statist. Assoc.}, 62:30--44, 1967.

\bibitem{Met10}
A.~Metzler.
\newblock On the first passage problem for correlated {B}rownian motion.
\newblock {\em Statistics and Probability Letters}, 80:277--284, 2010.

\bibitem{Pit96}
V.~I. Piterbarg.
\newblock {\em Asymptotic methods in the theory of {G}aussian processes and
  fields}, volume 148 of {\em Translations of Mathematical Monographs}.
\newblock American Mathematical Society, Providence, RI, 1996.

\bibitem{Pit15}
V.~I. Piterbarg.
\newblock {\em Twenty lectures about {G}aussian processes}.
\newblock Atlantic Financial Press, London, New York, 2015.

\bibitem{MR3493177}
Vladimir~I. Piterbarg.
\newblock High extrema of {G}aussian chaos processes.
\newblock {\em Extremes}, 19(2):253--272, 2016.

\bibitem{Res87}
S.~Resnick.
\newblock {\em Extreme {V}alues, {R}egular {V}ariation and {P}oint
  {P}rocesses.}
\newblock Springer-Verlag, 1987.

\bibitem{RS06}
L.~C.~G. Rogers and L.~Shepp.
\newblock The correlation of the maxima of correlated {B}rownian motions.
\newblock {\em J. Appl. Prob.}, 43(2):880--883, 2006.

\bibitem{ShaoWang13}
J.~Shao and X.~Wang.
\newblock Estimates of the exit probability for two correlated {B}rownian
  motions.
\newblock {\em Adv. Appl. Prob.}, 2013(45):37--50.

\bibitem{TeuGoo1994}
M.~Teunen and M.~Goovaerts.
\newblock Double boundary crossing result for the {B}rowian motion.
\newblock {\em Scandinavian Actuarial Journal}, 1994(2):139--150.

\bibitem{HH19}
R.~van~der Hofstad and H.~Honnappa.
\newblock Large deviations of bivariate {G}aussian extrema.
\newblock {\em Queueing Systems}, 93:333--349, 2019.

\bibitem{ZhouXiao17}
Y.~Zhou and Y.~Xiao.
\newblock Tail asymptotics for the extremes of bivariate {G}aussian random
  fields.
\newblock {\em Bernoulli}, 2017(23):1566--1598.

\end{thebibliography}
\newpage

\end{document}